\documentclass[12pt]{amsart}
\font\emailfont=cmtt10

\usepackage{amsmath,amsthm,amsfonts,amscd,flafter,epsf}

\newtheorem{theorem}{Theorem}[section]
\newtheorem{prop}[theorem]{Proposition}
\newtheorem{cor}[theorem]{Corollary}
\newtheorem{conj}[theorem]{Conjecture}
\newtheorem{lemma}[theorem]{Lemma}

\newtheorem{defn}[theorem]{Definition}

\def\endproof{\relax\ifmmode\expandafter\endproofmath\else
  \unskip\nobreak\hfil\penalty50\hskip.75em\hbox{}\nobreak\hfil\bull
  {\parfillskip=0pt \finalhyphendemerits=0 \bigbreak}\fi}
\def\endproofmath$${\eqno\bull$$\bigbreak}
\def\bull{\vbox{\hrule\hbox{\vrule\kern3pt\vbox{\kern6pt}\kern3pt\vrule}\hrule}}

\newcounter{bean}
\newcommand{\Q}{\mathbb{Q}}
\newcommand{\R}{\mathbb{R}}

\newcommand{\Z}{\mathbb{Z}}

\newcommand{\OneHalf}{\frac{1}{2}}

\newcommand{\Zmod}[1]{\Z/{#1}\Z}

\newcommand{\grad}{\vec\nabla}

\newcommand{\cm}{\cdot}

\newcommand{\Nbd}[1]{{\mathrm{nd}}(#1)}
\newcommand{\nbd}[1]{\Nbd{#1}}

\newcommand{\ModSWfour}{\mathcal{M}}
\newcommand{\ModFlow}{\ModSWfour}

\newcommand{\SpinC}{{\mathrm{Spin}}^c}

\newcommand\abuts\Rightarrow
\newcommand\Sym{\mathrm{Sym}}

\newcommand\RelSpinC{\underline{\SpinC}}
\newcommand\relspinc{\underline{\spinc}}

\newcommand\ts{\mathbf t}
\newcommand\x{\mathbf x}
\newcommand\w{\mathbf w}
\newcommand\z{\mathbf z}
\newcommand\p{\mathbf p}
\newcommand\q{\mathbf q}
\newcommand\n{\mathbf n}
\newcommand\y{\mathbf y}
\newcommand\m{\mathbf m}

\newcommand\ModSphere{\ModFlow\left({\mathbb S}\longrightarrow 
\Sym^{g-1}(\Sigma_{1})\times \Sym^2(\Sigma_{2})\right)}
\newcommand\ModSpheres\ModSphere

\newcommand\gr{\mathrm{gr}}
\newcommand\Mas{\mu}
\newcommand\UnparModSp{\widehat \ModSp}
\newcommand\UnparModFlow\UnparModSp
\newcommand\Mod\ModSp

\newcommand\PD{\mathrm{PD}}

\newcommand\Lk{\lk}

\newcommand{\spinc}{\mathfrak s}

\newcommand\ModMaps{\mathcal M}
\newcommand\ModSp\ModMaps

\newcommand\Ta{{\mathbb T}_{\alpha}}
\newcommand\Tb{{\mathbb T}_{\beta}}
\newcommand\Tc{{\mathbb T}_{\gamma}}

\newcommand\alphas{\mbox{\boldmath$\alpha$}}

\newcommand\betas{\mbox{\boldmath$\beta$}}
\newcommand\gammas{\mbox{\boldmath$\gamma$}}

\newcommand\Field{\mathbb F}

\newcommand\Dual{\mathcal D}
\newcommand\Duality\Dual

\newcommand\Gen{\mathfrak{S}}
\newcommand\HGrading{\mathfrak h}
\newcommand\U{{\mathbf 1}^*}
\newcommand\HH{\mathbb{H}}
\newcommand\lk{\mathrm{lk}}

\newcommand\orL{\vec{L}}

\newcommand\CFLa{\widehat{\mathrm CFL}}

\newcommand\HFLa{\widehat {\mathrm{HFL}}}

\newcommand\ws{\mathbf w}
\newcommand\zs{\mathbf z}
\newcommand\sss{\mathbf s}

\newcommand\spincrel\relspinc

\newcommand\CFK{\mathrm{CFK}}
\newcommand\HFK{\mathrm{HFK}}

\newcommand\CFKa{\widehat\CFK}

\newcommand\HFKa{\widehat\HFK}

\title[{Link Floer homology and the Thurston norm}] 
{Link Floer homology and the Thurston norm}

\author[Peter Ozsv{\'a}th]{Peter Ozsv\'ath}
\address{Department of
Mathematics, Columbia University, 
New York, NY 10027 \newline
\indent{\emailfont{petero@math.columbia.edu}}}
\thanks{PSO was supported by NSF grant number DMS-050581}

\author[Zolt{\'a}n Szab{\'o}]{Zolt{\'a}n Szab{\'o}} 
\address{Department of
Mathematics, Princeton University, New Jersey 08544 \newline
\indent{\emailfont{szabo@math.princeton.edu}}}
\thanks{ZSz was supported by NSF grant number DMS-0406155}

\begin{document}

\begin{abstract}  
  We show that link Floer homology detects the Thurston norm of a link
  complement. As an application, we show that the Thurston polytope of
  an alternating link is dual to the Newton polytope of its
  multi-variable Alexander polynomial.  To illustrate these
  techniques, we also compute the Thurston polytopes of several
  specific link complements.
\end{abstract}

\maketitle
\section{Introduction}

Heegaard Floer homology is an invariant of closed, oriented
three-manifolds which is defined using Heegaard diagrams of the
three-manifold~\cite{HolDisk}. The construction uses a suitable
variant of Lagrangian Floer homology in a symmetric product of a
Heegaard surface.  In~\cite{Knots} and~\cite{RasmussenThesis}, this
construction is refined to define knot Floer homology, an invariant
for null-homologous knots in an arbitrary (closed, oriented)
three-manifold. For the case of knots in the three-sphere, this
invariant is a bigraded Abelian group, whose graded Euler
characteristic is the Alexander polynomial.  Moreover, in this case,
knot Floer homology detects the genus of the knot~\cite{GenusBounds}.

In~\cite{Links}, the constructions from knot Floer homology are
generalized to the case of links in $S^3$. For an $\ell$-component
link, this gives a multi-graded Abelian group, with one grading for
each component of the link, and an additional grading (called the {\em
Maslov grading}). More precisely, let $\orL\subset S^3$ be an oriented
link, let $\mu_i$ be a meridian for the $i^{th}$ component $L_i$ of
$L$, and let $\HH\subset H_1(S^3-L;\R)$ be the affine lattice over
$H_1(S^3-L;\Z)$ given by elements $$\sum_{i=1}^\ell a_i\cm [\mu_i],$$
where $a_i\in\Q$ satisfies the property that $$2a_i+\lk(L_i,L-L_i)$$
is an even integer. Then, we have a finitely generated vector space
over $\Field=\Zmod{2}$ which splits as follows
$$\HFLa(\orL)=\bigoplus_{s\in \HH, d\in\Z} \HFLa_d(\orL,s).$$
The rank
of $\HFLa(\orL,h)$ is independent of the orientation of $L$ (cf.
Lemma~\ref{lemma:IndepOrientation} below), and hence the orientation
is often dropped from the notation.  The relationship with the
multi-variable Alexander polynomial $\Delta_L(T_1,...,T_\ell)$ is
given by the formula
\begin{equation}
  \label{eq:EulerChar}
\sum_{i=1}^\ell \chi(\HFLa_*(\orL,s))\cm e^s 
= \left(\prod_{i=1}^\ell (T_i^{\OneHalf}-T_i^{-\OneHalf})\right)\cm
\Delta_{L}(T_1,...,T_\ell),
\end{equation}
(see~\cite[Equation~\eqref{Links:eq:EulerHFLa}]{Links})
where here $s\mapsto e^s$ denotes the map from $\HH$
to Laurent polynomials 
$\Z[T_1^{\pm 1/2},...,T_{\ell}^{\pm 1/2}]$
which associates to the homology class
$$s=\sum_{i=1}^\ell a_i\cm [\mu_i]$$
the Laurent polynomial
$$T_1^{a_1}\cm...\cm T_\ell^{a_\ell}.$$

Our aim here is to extract topological information from these groups,
concerning the minimal genus of embedded surfaces representing a given
homology class. This information is neatly encoded in Thurston's
semi-norm on homology, cf.~\cite{ThurstonNorm}.

Recall that if $F$ is a compact, oriented, but possibly disconnected
surface-with-boundary $F=\bigcup_{i=1}^n F_i$, its {\em complexity}
is given by
$$\chi_-(F)=\sum_{\{F_i\big| \chi(F_i)\leq 0\}} -\chi(F_i).$$
Given any homology class $h\in H_2(S^3,L)$, it is easy to see that
there is a compact, oriented surface-with-boundary
embedded in $S^3-\nbd{K}$ representing $h$.  
Consider the function from $H_2(S^3,L;\Z)$ to the integers defined by
$$x(h)=\min_{\{F\hookrightarrow S^3-\nbd{K}\big|[F]=h\}} \chi_-(F).$$
According to Thurston~\cite{ThurstonNorm},
this can be naturally extended to a semi-norm, the {\em Thurston
semi-norm}, 
$$x\colon H_2(S^3,L;\R)\longrightarrow \R.$$

Link Floer homology also provides a function
$$y\colon H^1(S^3-L;\R) \longrightarrow \R$$ defined by the formula
$$y(h) =
\max_{\{s\in \HH\subset H_1(L;\R) \big| \HFLa(L,s)\neq 0\}} |\langle s, h\rangle|.$$

A link is said to have {\em trivial components} if it has some unknotted
component which is also unlinked from the rest of the link. Clearly,
adding trivial components does not change the Thurston semi-norm.

\begin{theorem}
  \label{thm:ThurstonNorm}
  The link Floer homology groups of an oriented link $\orL$ with no
  trivial components determines the Thurston norm of its complement, in the sense that
  for each $h\in H^1(S^3-L;\R)$
  $$x(\PD[h])
  +\sum_{i=1}^\ell |\langle h, \mu_i\rangle |
  =2y(h),$$
  where $\mu_i$ is the meridian for the $i^{th}$ component of $L$,
  so that $|\langle h, \mu_i\rangle|$ denotes the absolute value of the Kronecker
  pairing of $h\in H^1(S^3-L;\R)$ with the homology class $\mu_i$, thought of as an element
  of $H_1(S^3-L;\R)$.
\end{theorem}

Let $V$ be a finite-dimensional vector space equipped with a
(semi-)norm $N$ which is linear on rays in $V$. Such a (semi-)norm is
determined by its unit ball which, in the case of $x$, $y$, and
$|\cdot |$ are polytopes. Moreover, it is sometimes
useful to think about the dual norm $N^*$,
$$N^*(\xi)=\sup_{\{v\in V\big|N(v)=1\}}|\xi(v)|.$$ The unit ball of $N^*$ is
the dual of the unit ball of $N$ (in particular, the faces of one 
correspond to the vertices of the other).

The unit ball for $x^*$ is a polytope in $H_1(S^3-L;\R)$, called the
{\em dual Thurston polytope}. For $y$, we obtain the {\em link Floer
homology polytope}, which is the convex hull of those $s\in\HH$ for
which $\HFLa(L,s)\neq 0$.  Theorem~\ref{thm:ThurstonNorm} says, then,
that twice the link Floer homology polytope is the set of points which
can be written as a sum of an element of the dual Thurston polytope
and an element of the symmetric hypercube in $H^1(S^3-L)$ with
edge-length two.

Theorem~\ref{thm:ThurstonNorm} has a number of antecedents. Monopole
Floer homology~\cite{KMbook} detects the Thurston norm of a closed
three-manifold, according to a fundamental result of Kronheimer and
Mrowka~\cite{KMthurston}, see also~\cite{KMOSz}, building on results
of Gabai~\cite{DiskDecomposition} and
Eliashberg-Thurston~\cite{EliashbergThurston}.  In a similar manner,
Heegaard Floer homology, and also Floer homology for knots, detects
the corresponding Thurston norms according to~\cite{GenusBounds},
building on further results in topology and symplectic geometry,
notably~\cite{GabaiKnots}, \cite{Eliashberg}, \cite{Etnyre},
\cite{DonaldsonLefschetz}.  A generalization of this result to links
has been established by Ni~\cite{Ni}.  His theorem amounts to
Theorem~\ref{thm:ThurstonNorm} for the case of where $h$ is one of the
$2^\ell$ cohomology classes with $|\langle h,\mu_i\rangle| = 1$ for
$i=1,...,\ell$.  In fact, our proof of Theorem~\ref{thm:ThurstonNorm}
reduces to this case, in view of properties of both $x$ and $y$ under
cabling, compare also~\cite{EisenbudNeumann} and~\cite{Hedden}
respectively.

In a slightly different direction, it is a classical fact that the
degree of the Alexander polynomial gives a lower bound on the genus of
a knot. In~\cite{McMullen}, McMullen generalizes this result, showing
that the Newton polytope of the multi-variable Alexander polynomial is
contained in the dual Thurston polytope. 

At present, there is no algorithm for calculating link Floer homology
in general. However, there are some useful calculational devices, such
as skein exact sequences, and also in some cases, link Floer homology
can be calculated directly by examining Heegaard diagrams.  In
practice it is typically much easier to calculate the link Floer
homology polytope than the full link Floer homology.

We have the following result for
alternating links:

\begin{theorem}
  \label{thm:AltLinks}
  Let $L$ be a link with connected, alternating projection.  The rank
  of $\HFLa(L,s)$ is the absolute value of the coefficient of $e^s$ in
  $$\left(\prod_{i=1}^\ell (T_i^{\OneHalf}-T^{-\OneHalf})\right)\cm
\Delta_{L}(T_1,...,T_\ell).$$
\end{theorem}

A more precise version is stated
in~\cite[Theorem~\ref{Links:thm:AltLink}]{Links}, which in turn
follows rather quickly from results of~\cite{AltKnots}.

Combining Theorems~\ref{thm:ThurstonNorm} and~\ref{thm:AltLinks}, we obtain the
following generalization of a classical theorem of Crowell and
Murasugi~\cite{Crowell}, \cite{Murasugi}, affirming a conjecture of
McMullen~\cite{McMullen}:

\begin{cor}
  \label{cor:AltLinks} Let $L\subset S^3$ be a link with $\ell$
  components which admits a connected, alternating
  projection. Consider the convex hull of all the points in
  $\ell$-dimensional space which correspond to non-zero terms in the
  multi-variable Alexander polynomial of $L$ (i.e. the Newton polytope
  of the multi-variable Alexander polynomial). This polytope, scaled
  by a factor of two, is the dual Thurston polytope of the complement
  of $L$.
\end{cor}

As a further illustration, we also calculate the Thurston polytopes of
various links. Specifically, we describe the Thurston polytopes of
those nine-crossing links which were not described in
~\cite{McMullen}, namely, $9^2_{41}$ (which is alternating),
$9^2_{50}$, and $9^3_{15}$. Then, we turn our attention to a
two-component link with trivial Alexander polynomial, the
$10$-crossing Kinoshita-Terasaka link.

Of course, the present paper depends on the link Floer homology
of~\cite{Links}. With this said, it is worth underscoring the fact
that we use here only a very minimal version of link Floer homology:
the more complicated gluing results for pseudo-holomorphic curves are
not needed in our applications.

In~\cite{ThurstonNorm}, Thurston shows that the set of elements of
$H^1(S^3-L;\Z)$ which represent fibrations of the link complement
correspond to certain open faces of his polytope, called  {\em fibered faces}.
According to Theorem~\ref{thm:ThurstonNorm}, to each extremal point
$P=\sum_{i=1}^\ell a_i\cm \mu_i \in H_1(S^3-L;\Z)$ of the dual
Thurston polytope, there is a set $s(P)$ of corresponding extremal points in
the link Floer homology polytope; these are the extremal points which can be written as
$(P+\sum_{i=1}^\ell \epsilon_i \cm \mu_i^*)/2$ (where here $\{\mu_i^*\}_{i=1}^\ell$
denotes the dual basis in $H^1(S^3-L;\R)$ for $\{\mu_i\}_{i=1}^\ell H_1(S^3-L;\R)$).
The methods in the proof of Theorem~\ref{thm:ThurstonNorm} readily
give the following simple geometric consequence for these groups:

\begin{prop}
        \label{prop:Fibrations}
        If $P\in H^1(S^3-L;\Z)$ corresponds to a fibered face of the
        Thurston polytope, then for each $h\in s(P)$, $\HFLa(S^3,h)$
        is one-dimensional.
\end{prop}

Conversely, one is inclined to believe the following:

\begin{conj}
  \label{conj:Fibrations}
  If $P$ corresponds to a face of the Thurston polytope with the
  property that for some $h\in s(P)$, $\HFLa(S^3,h)$ is
  one-dimensional, then $P$ corresponds to a fibered face.
\end{conj}

An analogous conjecture has been made for knots~\cite{Survey}. 

In Section~\ref{sec:Background}, we give some of the background for
the link Floer homology from~\cite{Links}, with a special emphasis on
the part of the theory relevant to us for our present purposes. In
Section~\ref{sec:Proof}, we prove Theorem~\ref{thm:ThurstonNorm} and
Corollary~\ref{cor:AltLinks}. In Section~\ref{sec:Examples}, we turn
to some applications, and some illustrative calculations (involving
links with trivial Alexander polynomial).  We conclude with a proof of
Proposition~\ref{prop:Fibrations} in Section~\ref{sec:Fibrations}.

\subsection{Acknowledgements}
We would like to thank David Gabai, Matthew Hedden, Walter Neumann, Yi
Ni, and Jacob Rasmussen for interesting conversations during the
course of this work. We are particularly indebted to Jake for his many
valuable suggestions following a thorough reading an early version of
this paper.

\section{Background on Link Floer homology}
\label{sec:Background}

\subsection{Definitions}
Link Floer homology is defined in a fairly general context
in~\cite{Links} (compare also~\cite{Knots} and~\cite{RasmussenThesis}
for the case of knots). We sketch here the parts of this construction
which we need presently.

Given an oriented surface $\Sigma$ of genus $g$ and a positive integer
$\ell$, a $g+\ell-1$-tuple of embedded, disjoint curves whose homology
classes span a $g$-dimensional subspace of $H_1(\Sigma)$ specifies a
handlebody which is bounded by $\Sigma$.  Fix, then, two such
$g+\ell-1$-tuples of circles
$\alphas=\{\alpha_1,...,\alpha_{g+\ell-1}\}$ and
$\betas=\{\beta_1,...,\beta_{g+\ell-1}\}$, and let $U_\alpha$ and
$U_\beta$ denote the corresponding handlebodies. Fix also
$2\ell$-points in
$$\Sigma-\alpha_1-...-\alpha_{g+\ell-1}-\beta_1-...-\beta_{g+\ell-1},$$
denoted $\ws=\{w_1,...,w_\ell\}$ and $\zs=\{z_1,...,z_\ell\}$.
Suppose that $w_i$ and $z_i$ can be connected by arcs $$\xi_i\subset
\Sigma-\alpha_1-...-\alpha_{g+\ell-1}$$ and $$\eta_i\subset
\Sigma-\beta_1-...-\beta_{g+\ell-1}.$$

In this case, we can specify a link in $Y=U_\alpha\cup_{\Sigma}
U_\beta$ as follows.  Let $\xi_i'$ denote the arc in $U_\alpha$
obtained by pushing $\xi_i$ into the handlebody so that it meets
$\Sigma$ only at its endpoints $w_i$ and $z_i$, and let $\eta_i'$
denote the analogous push-off of $\eta_i$ in $U_\beta$.
Our link $L\subset Y$, then is given by
$$\bigcup_{i=1}^{\ell} \xi_i'\cup \eta_i'.$$
An orientation for $Y$ is
inherited from the orientation of $U_\alpha\subset Y$, which in turn
is oriented so that the given orientation on $\Sigma$ agrees with the
orientation it inherits from being the boundary of $U_\alpha$.
Moreover, an orientation for $\orL$ is specified by the convention
that the subarc $\xi_i'\subset \orL$ inherits an orientation as a path
from $w_i$ to $z_i$.

In this case, we say that $(\Sigma,\alphas,\betas,\ws,\zs)$ is
a {\em $2\ell$-pointed Heegaard diagram compatible with the
oriented link $\orL\subset Y$}.

For our applications, we restrict attention to the case where $Y\cong
S^3$.

A {\em periodic domain} for a $2\ell$-pointed Heegaard diagram is a
sum $P=\sum {m_i} D_i$, where here $D_i$ are the closures of the
components of
$$\Sigma-\alpha_1-...-\alpha_{g+\ell-1}-\beta_1-...-\beta_{g+\ell-1},$$
and with the additional properties that 
$$\partial P = \sum a_i\cm \alpha_i + \sum b_i\cm \beta_i,$$
and whose local multiplicities at each of the $w_i$ and $z_i$ are zero.
A $2\ell$-pointed Heegaard diagram is called {\em admissible}
if each non-zero periodic domain $P$ has at least one positive
and at least one negative local multiplicity  ($m_i$). 

Given a $2\ell$-pointed Heegaard diagram, we can form the
$g+\ell-1$-fold symmetric product of the Heegaard surface
$\Sym^{g+\ell-1}(\Sigma)$, equipped with the pair of tori
\begin{eqnarray*}
\Ta=\alpha_1\times...\times\alpha_{g+\ell-1}
&{\text{and}}&
\Tb=\beta_1\times...\times\beta_{g+\ell-1}.
\end{eqnarray*}
Let $\Gen$ denote the set of intersection points between $\Ta$ and $\Tb$.

Link Floer homology~\cite{Links} is a version of Lagrangian Floer
homology~\cite{FloerLag}, \cite{FloerUnregularized} in this context.
Specifically, starting from an admissible $2\ell$-pointed Heegaard
diagram for a link, where all the curves $\alpha_i$ and $\beta_j$ meet
transversally, we consider the chain complex $\CFLa$ generated as a
vector space over $\Field$ by the intersection points $\Gen$,
endowed with the differential
\begin{equation}
\label{eq:DefD}
\partial\x = \sum_{\y\in\Gen}\sum_{\{\phi\in\pi_2(\x,\y)\big|
  n_{\ws}(\phi)=n_{\zs}(\phi)=0,\Mas(\phi)=1\}}
\#\left(\frac{\ModFlow(\phi)}{\R}\right)\y.
\end{equation}
Here, $\pi_2(\x,\y)$ is the space of homology classes of Whitney disks
connecting $\x$ and $\y$, $n_{\ws}(\phi)\in\Z^\ell$ is the
$\ell$-tuple $(n_{w_1}(\phi),...,n_{w_\ell}(\phi))$, where
$n_{w_i}(\phi)$ denotes the algebraic intersection number of $\phi$
with $\{w_i\}\times \Sym^{g+\ell-2}(\Sigma)\subset
\Sym^{g+\ell-1}(\Sigma)$, $n_{\z}(\phi)$ is defined analogously,
$\ModFlow(\phi)$ denotes the moduli space of pseudo-holomorphic
representatives of $\phi$, and $\Mas(\phi)$ denotes its expected
dimension. The quantity $\#(\frac{\ModFlow(\phi)}{\R})$ denotes the
number of points in this finite set, counted modulo two. When the
pseudo-holomorphic condition is suitably generic, we have that
$\partial ^2 = 0$, i.e. $\CFLa$ is in fact a chain complex.

In fact, the chain complex $\CFLa$ can be endowed with a relative
Maslov grading, specified by
$$\gr(\x)-\gr(\y)=\Mas(\phi)-2\sum_{i=1}^\ell n_{w_i}(\phi),$$ where
$\phi$ is any disk from $\x$ to $\y$. Note that, as the notation
suggests, this quantity is independent of the particular choice of
$\phi$.  With this convention, then, $\CFLa$ inherits a relative
$\Z$-grading, with the property that the boundary operatory of
Equation~\eqref{eq:DefD} drops grading by one. In fact, this relative
grading can be enhanced to an absolute $\Z$-grading (the Maslov
grading) as well, but we have no need for this additional structure in
the present paper.

We can define a function
$$\HGrading_{\w,\z}\colon \Gen \longrightarrow \HH$$
with the property that
\begin{equation}
\label{eq:DefHGrading}
\HGrading_{\w,\z}(\x)-\HGrading_{\w,\z}(\y)=\sum_{i=1}^\ell
(n_{z_i}(\phi)-n_{w_i}(\phi))\mu_i,
\end{equation}
where $\phi\in\pi_2(\x,\y)$ is
any Whitney disk connecting $\x$ and $\y$.  We have a splitting of
$\CFLa$ into summands indexed by homology classes $h\in \HH$,
generated by those intersection points $\x$ with 
$\HGrading_{\w,\z}(\x)=h$. 

The homology group of this summand is the {\em link Floer homology group}
of $L$,
$\HFLa(L,h)$; we can collect these into one group by
$$\HFLa(L)=\bigoplus_{h\in \HH} \HFLa(L,h).$$

As the notation suggests, this is a link invariant, according to one
of the main results of~\cite{Links}.

Strictly speaking, the function $\HGrading_{\w,\z}$ is characterized by
Equation~\eqref{eq:DefHGrading} only up to an overall translation.
We describe how to remove this ambiguity with the help
of a symmetry, cf. Equation~\eqref{eq:Symmetry}. An alternative
approach proceeds via the notion of ``relative $\SpinC$ structures'',
which we recall in Subsection~\ref{subsec:RelSpinC}.

\subsection{Symmetries}

Heegaard Floer homology for links enjoys a number of basic properties.
For example, its Euler characteristic is determined by the multi-variable
Alexander polynomial, as in Equation~\eqref{eq:EulerChar}.
Another fundamental property is the following isomorphism
of {\em relatively graded} $\Z$-graded groups
(generalizing the usual symmetry of the Alexander polynomial):
\begin{equation}
  \label{eq:Symmetry}
  \HFLa_*(\orL,h)\cong \HFLa_*(\orL,-h),
\end{equation}
which holds for any fixed $h\in \HH$, 
see~\cite[Equation~\eqref{Links:eq:SymmetryH}]{Links}.

\begin{lemma}
        \label{lemma:IndepOrientation}
        Let $\orL_1$ and $\orL_2$ denote two different orientations
        on the same underlying link $L$. Then, for each $h\in \HH$,
        there is an isomorphism of relatively $\Z$-graded groups
        $$\HFLa_*(\orL_1,h)\cong \HFLa_*(\orL_2,h).$$
\end{lemma}

\begin{proof}
        Consider a $2\ell$-pointed Heegaard diagram for the oriented
        link $\orL_1$. Given any other orientation $\orL_2$ on the
        same underlying link, we can obtain a corresponding
        $2\ell$-pointed Heegaard diagram for $\orL_2$ by reversing the
        roles of some pairs of the $w_i$ and $z_i$. Obviously the
        differential in Equation~\eqref{eq:DefD} is unchanged by this
        operation. Thus, the total rank of $\HFLa$ is independent
        of the orientation used on the link. 

        Next, we consider the splitting of this group into components
        indexed by elements of $\HH$.  Letting \begin{eqnarray*}
        \HGrading_1\colon \Gen\longrightarrow \HH &{\text{and}}&
        \HGrading_2\colon \Gen\longrightarrow \HH
        \end{eqnarray*} be the maps for these two choices of $\ws$ and
        $\zs$, we see that
        $$\HGrading_1(\x)-\HGrading_1(\y)=\HGrading_2(\x)-\HGrading_2(\y)$$
        for any $\x, \y\in\Gen$.  This follows at once from
        Equation~\eqref{eq:DefHGrading}: we use one homology class
        $\phi\in\pi_2(\x,\y)$ to calculate both sides, and observe
        that reversing the orientation of the $i^{th}$ component
        changes at once the sign of $n_{z_i}(\phi)-n_{w_i}(\phi)$ and
        also the sign of the $i^{th}$ meridian $\mu_i$. It follows
        that there is some fixed $h\in H_1(S^3-L;\Q)$ with the
        property that for all $\x\in\Gen$,
        $\HGrading_1(\x)=\HGrading_2(\x)+h$. By symmetry
        (Equation~\eqref{eq:Symmetry}), it follows that $h=0$, and the
        lemma is complete.
\end{proof}

The {\em absolute} $\Z$-grading on $\HFLa_*$ does, however, depend on
the orientation of $L$. But the Floer homology polytope depends only
on the set of $h$ with non-trivial  $\HFLa_*(\orL,h)$ which, according to
Lemma~\ref{lemma:IndepOrientation} is independent of the orientation
on $L$. Indeed, we will think of link Floer homology only with its
relative Maslov grading, and hence we will often drop the orientation
of $L$ from the notation.

\subsection{Relationship with knot Floer homology}

The construction of Heegaard Floer homology for knots predates the
corresponding construction for links~\cite{Knots},
\cite{RasmussenThesis}. Moreover, the paper~\cite{GenusBounds}, contains
a proof of Theorem~\ref{thm:ThurstonNorm} for the case of knots.
Specifically, it is shown there that if $K$ is a knot, then the
minimal genus of any Seifert surface for the knot, its {\em Seifert
  genus} $g(K)$, is given by
$$\max_{\{s\in \Z\big| \HFKa(K,s)\neq 0\}} |s|.$$
For the case of knots, we write $\HFKa$ for the corresponding link
Floer homology, which we think of as $\Z$-graded, under some
identification $\Z\cong H_1(S^3-K;\Z)$.

In fact, as described in
\cite[Proposition~\ref{Knots:prop:LinksToKnots}]{Knots}, since the
knot invariant can be defined for null-homologous knots in an
arbitrary three-manifold, it can also be used to define an invariant
for oriented links in $S^3$, in the following manner.  Starting form
an oriented link $\orL$ in $S^3$ with $\ell$ components, we attach
$\ell-1$ one-handles to $S^3$, simultaneously attaching one-handles to
our link, so as to obtain a connected knot $\kappa(\orL)$ inside
$\#^{\ell-1}(S^2\times S^1)$.  We then define the ``knot Floer
homology'' for the oriented link $\orL\subset S^3$ to be the Floer
homology of the associated knot $\kappa(\orL)\subset
\#^{\ell-1}(S^3\times S^1)$, written
$$\HFKa(\orL)=\bigoplus_{s\in\Z}\HFKa(\orL,s).$$ Note that the graded
Euler characteristic of this theory is a (suitably normalized) version
of the Alexander-Conway polynomial,
cf.~\cite[Equation~\eqref{Knots:eq:EulerChar}]{Knots}.
 
In~\cite{Ni}, Ni shows that the breadth of these homology groups calculates the
Seifert genus of the oriented link, in the following sense.

\begin{theorem} (Ni~\cite{Ni})
\label{thm:Ni}
        Fix an oriented link $\orL$ with $\ell$ components. Then,
        $$2\max{\{s\in\Z\big|\HFKa(\orL,s)\neq 0\}}
        =\min_{\{F\hookrightarrow S^3\big|\partial F=\orL\}} \ell-\chi(F).$$
\end{theorem}

The knot Floer homology of $\kappa(\orL)$ and the link Floer homology
of $\orL$ can be immediately related by the following:

\begin{lemma}
        \label{lemma:HFKandHFL}
  There is a spectral sequence whose $E_2$ term is
  $$\sum_{a_1+...+a_\ell=s} \HFLa(\orL,\sum_{i=1}^\ell a_i\cm \mu_i)$$
  and whose $E^{\infty}$ term is $\HFKa(\orL,s)$.
\end{lemma}

\begin{proof}
  Start from a pointed Heegaard diagram for $\orL$. By attaching
  one-handles to the surface, connecting $z_i$ to $w_{i+1}$ for
  $i=1,...,\ell-1$, and forgetting all the basepoints except $w_1$ and
  $z_\ell$, we obtain a doubly-pointed Heegaard diagram for
  $\kappa(\orL)\subset \#^{\ell-1}(S^2\times S^1)$. The remaining basepoints
  $z_1,...,z_{\ell-1}$ can be thought of as giving a further $\Z^{\ell-1}$
  filtration of the chain complex $\CFKa(\#^{\ell-1}(S^2\times S^1),\kappa(\orL),s)$.
  The associated graded object for this filtration is 
  $$\bigoplus_{\{(a_1,...,a_{\ell})\big|\sum_{i=1}^\ell a_i = s\}} \HFLa(\orL,
  \sum_{i=1}^\ell a_i\cm \mu_i).$$
  The lemma now follows from the Leray spectral sequence of this filtration.
\end{proof}

In fact, in Theorem~\ref{Links:thm:IdentifyWithLinkHomology}
of~\cite{Links} more is proved: it is shown that the above spectral
sequence collapses, so that
$$
\HFKa(\orL,s)\cong \bigoplus_{a_1+...+a_\ell=s} \HFLa(\orL,\sum_{i=1}^\ell a_i\cm \mu_i).$$
We will not need this stronger form in the present applications;
Lemma~\ref{lemma:HFKandHFL} suffices. Indeed, it will be useful
to have the following combination of the lemma with Ni's theorem:

\begin{prop}
\label{prop:NiPlusSS}
Let $\orL$ be an oriented link, and let 
$$m=
\max_{\{h=\sum_{i=1}^\ell a_i\cm \mu_i\in H_1(S^3-L)\big|\HFLa(L,h)\neq 0\}} 
 \sum a_i.$$
Suppose moreover that 
there is a unique $h=\sum_{i=1}^\ell a_i\cm \mu_i\in H_1(S^3-L)$ with $\HFLa(S^3-L,h)\neq 0$ for which 
$\sum_{i=1}^\ell a_i = m$.
Then, 
$$2m=\min_{\{F\hookrightarrow S^3\big|\partial F=\orL\}} \ell-\chi(F).$$
\end{prop}

\begin{proof}
  By our hypothesis, the $E_2$ term in the spectral
  sequence from Lemma~\ref{lemma:HFKandHFL} 
  converging to $\HFKa(\orL,m)$ consists of the
  single term $\HFLa(\orL,h)$, and hence it collapses;
  i.e. 
  $$\HFKa(\orL,m)\cong \HFLa(\orL,h).$$ Note also that for
  all $s>m$, the $E_2$ term of the spectral sequence converging
  to $\HFKa(\orL,s)$ vanishes. Thus, we have that 
  $m=\max{\{s\in\Z\big|\HFKa(\orL,s)\neq 0\}}$,
  and the lemma now follows from Ni's theorem.
\end{proof}

\subsection{Relative $\SpinC$ structures}
\label{subsec:RelSpinC}

There is a conceptually more satisfying, if less practical, method of
thinking about the $\HH$-grading on link Floer homology, which is to
employ the notion of {\em relative $\SpinC$ structures} on the link
complement (cf. Section~\ref{Links:sec:HeegaardDiagrams}
of~\cite{Links}). 

Let $(M,\partial M)$ be an oriented three-manifold whose boundary
consists of a union of tori $T_1\cup...\cup T_\ell$. On each torus,
there is a preferred isotopy class of nowhere vanishing vector field,
containing those which are invariant under translation on the torus.
Consider nowhere vanishing vector fields on $M$ whose restriction to
$\partial M$ are tangent to the boundary, where they are
translationally invariant.  Following Turaev~\cite{Turaev}, we say
that $v$ and $v'$ are {\em homologous} if there is a ball $B\subset
M-\partial M$ with the property that the restrictions of $v$ and $v'$
to $M-B$ are homotopic through nowhere vanishing vector fields which
are tangent to $\partial M$.  The set of homology classes of such
vector fields is called the set of {\em relative $\SpinC$ structures},
and it is an affine space for $H^2(M,\partial M;\Z)$.  We denote this
set by $\RelSpinC(M,\partial M)$. In the case where $M=S^3-\nbd{L}$,
we denote the set by $\RelSpinC(S^3,L)$.

There is a natural map $$c_1\colon \RelSpinC(M,\partial M)
\longrightarrow H^2(M,\partial M),$$
induced by taking the nowhere
vanishing vector field $v$ to the first Chern class of the orthogonal
complement of $v$, relative to the natural trivialization on the
boundary given by outward pointing vectors.

Let $(\Sigma,\alphas,\betas,\ws,\zs)$ be a $2\ell$-pointed Heegaard
diagram for an oriented link $\orL$.  Given an intersection point
$\x\in\Gen$, we can define the associated relative $\SpinC$
structure $\relspinc_{\w,\z}(\x)\in \RelSpinC(S^3,L)$ as follows. Let
$f\colon S^3\longrightarrow \R$ be a self-indexing Morse function and
$g$ be a Riemannian metric on $S^3$ with the following properties:
\begin{itemize}
\item $f$ has $\ell$ index zero and index three
critical points, and $g+\ell$ index one and two cricial points, and mid-level $\Sigma$,
\item  $\alpha_i$ is the set of points flowing
into $\Sigma$ out of the $i^{th}$ index one critical point, and
$\beta_j$ is the set of points flowing into the $j^{th}$ index two
critical point,
\item the set of flowlines which pass through $\{w_i,z_i\}_{i=1}^\ell$
is identified with $L\subset S^3$, oriented so that 
$\orL$ is oriented upward at each $z_i$ (and downward at each $w_i$).
\end{itemize}
Such a Morse function is said to be {\em compatible} with the pointed
Heegaard diagram $(\Sigma,\alphas,\betas,\ws,\zs)$. Given
$\x\in\Gen$, we can consider the corresponding $g+\ell$-tuple of
gradient flowlines $\gamma_\x$ which connect the various index one and
two critical points, and the $\ell$-tuple of gradient flowlines
$\gamma_\w$ connecting the various index zero and three critical
points (and passing through all the $w_i$).  We can now construct a
nowhere vanishing vector field over $S^3$ with a closed orbit given by
$\orL$, by modifying the gradient vector field of $f$ in a
sufficiently small neighborhood of $\gamma_\x\cup \gamma_\w$.
The modification made in a neighborhood of $\gamma_\w$ is
concretely specified in Figure~\ref{Links:fig:VectorField} of~\cite{Links}
(and is not of primary importance to us at present).

Such a vector field can be viewed as a vector field on $S^3-\nbd{L}$
which is tangent to the boundary. The homology class of this vector
field induces a well-defined map
$$\relspinc_{\w,\z}\colon \Gen\longrightarrow \RelSpinC(S^3,L).$$

The relationship between this map and the map $\HGrading_{\w,\z}$ is given
by the formula
\begin{equation}
  \label{eq:RelSpinCFormula}
  c_1(\relspinc_{\w,\z}(\x))+\sum_{i=1}^{\ell}\PD[\mu_i] 
  = 2\cm \PD[\HGrading_{\w,\z}(\x)]
\end{equation}
where here we are using the Poincar\'e duality isomorphism
$$\PD\colon H_1(S^3-\nbd{L})\longrightarrow H^2(S^3,L).$$

\section{Proof of Theorem~\ref{thm:ThurstonNorm}.}
\label{sec:Proof}

To establish Theorem~\ref{thm:ThurstonNorm}, we compare how both $x$
and $y$ transform under cabling.  Before describing this, we introduce
some notation.

We have a basis for $H_1(\partial\nbd{L};\Z)$ given by
$\lambda_1,...,\lambda_\ell,\mu_1,...,\mu_\ell$, where $\lambda_i$ is
the longitude of the $i^{th}$ component of $L$, and $\mu_i$ is its
meridian. Correspondingly, given ${\bf p}=(p_1,...,p_{\ell})$ and
${\bf q}=(q_1,...,q_\ell)$, we can form a new link, the {\em cable}
$C_{\p,\q}(L)$ of $L$. This is the link gotten by inserting $\ell$
solid tori in $S^3-\nbd{L}$, where the $(p_i,q_i)$-torus knot (or
link) is contained in the solid torus inserted into the $i^{th}$
component of $\ell$. (Note that the number of components of the
$(p,q)$ torus link is $\gcd(p,q)$.)  An orientation on $L$ naturally
induces an orientation on the cable $C_{\p,\q}(L)$.

Given any ${\bf p}=(p_1,...,p_\ell)$, there is a unique choice
$Q(\p)=(Q_1,...,Q_\ell)$ with the property that
$$\sum_{i=1}^\ell Q_i\cm\mu_i + p_i\cm\lambda_i=0$$
as homology classes in $H_1(S^3-L)$;
specifically
\begin{equation}
  \label{eq:DefQi}
  Q_i=-\sum_{j\neq i} p_j\cm\lk(L_i,L_j).
\end{equation}

Let $j\colon S^3-\nbd{L} \longrightarrow S^3-C_{\p,\q}(L)$ denote the
natural inclusion map, and consider the induced maps
\begin{eqnarray*}
j_*\colon H_1(S^3-\nbd{L}) &\longrightarrow& H_1(S^3-C_{\p,\q}(L)) \\
j^*\colon  H^1(S^3-C_{\p,\q}(L)) &\longrightarrow &
H^1(S^3-\nbd{L}).
\end{eqnarray*}
In the case where each $q_i$ is relatively prime to $p_i$, the
components of $C_{\p,\q}(L)$ are in one-to-one correspondence with the
components of $L$. In this case, letting $\mu_i'$ be the meridian
of the $i^{th}$ component of $C_{\p,\q}(L)$, we clearly have that
$$j_*(\mu_i)=p_i\cm \mu_i'.$$

\begin{defn}
  \label{def:Minimal}
  Any one-dimensional homology class in the two-torus $T$ can be
  represented by an embedded, oriented one-manifold $C\subset T$. We
  say that the representative $C$ is {\em minimal} if no component is
  null-homologous, and any two of its components are orientation-preserving
  isotopic.
\end{defn}

It is well-known that the Thurston norm of $L$ can be understood in
terms of the minimal genus Seifert surfaces of its cables. For a
general discussion, see~\cite{EisenbudNeumann}. We recall this result
in the form we need it in the following:

\begin{lemma}
  \label{lemma:ThurstonNormCable} 
  Let $L$ be a link with no trivial components (i.e. no component of $L$
  is bounded by a disk which is disjoint from the rest of the link).
  Fix $\p=(p_1,...,p_\ell)$, where $p_i$ are positive integers, and
  let $Q_i$ be the corresponding integers as in
  Equation~\eqref{eq:DefQi}.  Then, for any $\ell$-tuple of integers
  $$\q=(q_1,...,q_\ell)$$
  with each $q_i\geq Q_i$, we have that
  $$x(C_{\p,\q}(L),\U)=x(L,j^*(\U))+\sum_{i=1}^{\ell}(q_i-Q_i)(p_i-1),$$
  where $\U\in H^1(S^3-C_{\p,\q}(L))$ denotes the cohomology class
  whose value on each oriented meridian for $C_{\p,\q}(L)$ is one.
\end{lemma}

\begin{proof}
  For $i=1,...,\ell$, let $T_i\subset S^3-C_{\p,\q}(L)$ be the torus
  which forms the boundary of a neighborhood of the $i^{th}$ component
  of $L$.  We claim that for each $\xi\in H^1(S^3-C_{\p,\q}(L))$,
  there is an embedded surface
  $$(F,\partial F)\hookrightarrow (S^3,C_{\p,\q}(L))$$
  of minimal
  complexity representing $\PD[\xi]$ with the property that $F$ meets
  each $T_i$ transversally and each intersection $T_i\cap F$ is
  minimal, in the sense of Definition~\ref{def:Minimal}.
  
  We arrange this as follows. Start from a minimal complexity surface
  $F_1$ meeting each $T_i$ transversally. Next, remove all the
  null-homotopic components of $F_1\cap T_i$, as follows. Suppose
  there is a circle $C_1\subset F_1\cap T_i$ which is null-homotopic
  in $T_i$.  Then, there is an innermost one $C_2$ (i.e. the disk in
  $T_i$ bounded by $C_2\subset F_1\cap T_i$ does not contain any other
  component of $F_1\cap T_i$). Surgering out this circle 
  gives a new embedded surface homologous to $F_1$ 
  whose complexity is no greater than that of $F_1$. We
  proceed in this manner until we obtain a complexity-minimizing
  surface $F_2$ for the homology class $\xi$ with the additional
  property that $F_2\cap T_i$ contains no null-homotopic components.
  
  Note now that $F_2$ is a complexity-minimizing surface representing
  $\PD[\xi]$ with the property that for each $i$, $F_2\cap T_i$
  consists parallel copies of the same (homotopically non-trivial)
  curve in $T_i$. Suppose next that there are two components $C_1$ and
  $C_2$ of $F_2\cap T_i$ which are oriented oppositely. We can then
  cut to obtain a new representative $F_3$ which meets $T_i$ in two
  fewer components. The Euler characteristic of $F_3$ agrees with that
  of $F_2$, and indeed its complexity must agree with that of $F_2$
  except in the special case where a sphere was created by the
  cut-and-paste operation. But this is possible only if $C_1$ and
  $C_2$ bound a disk on either side of $T_i$ in
  $(S^3-C_{\p,\q}(L))-T_i$. But this is impossible: $T_i$ is
  incompressible on both sides (we are using here the hypothesis
  that each $p_i$ is non-zero and that $L$ has no trivial components).
  
  Proceeding in this manner, we obtain a complexity-minimizing
  representative $F'$ for the homology class with the property that
  $T_i\cap F'$ is minimal. The $T_i$ divide $F'$ into surfaces
  $A$ in $S^3-\nbd{L}$ which represents $\PD[j^*(\xi)]$, and a
  collection of surfaces $B_i$ supported inside the solid tori bounded
  by $T_i$. The same arguments as above show that $A$ and $B_i$ are
  all complexity-minimizing in their respective relative homology
  classes.
  
  Specifically, $F'\cap T_i$ is the $(p_i,Q_i)$ torus link.  It is
  easy to see that the
  minimal complexity surface in the solid torus which meets its
  boundary in the $(p_i,Q_i)$ torus link, and whose other boundary
  component is the $(p_i,q_i)$ torus link inside has complexity
  $(q_i-Q_i)(p_i-1)$.
\end{proof}

In~\cite{Hedden}, Hedden studies the behaviour of knot Floer homology
under cabling. Among other things, he shows that the topmost
(non-trivial) Floer homology group of a sufficiently twisted cable of
a knot is isomorphic to the topmost knot Floer homology group of the
original knot. (See also~\cite{NiSutured} for a generalization of this
to other kinds of satellites.)  Adapting this to the context of link
Floer homology, we obtain the following:

\begin{prop}
  \label{prop:Cabling} Let $\p=(p_1,...,p_\ell)$ be an $\ell$-tuple
  of positive integers, each of which is greater than one.  Consider
  the cohomology class
  $$\theta=\sum_{i=1}^\ell p_i\cm \mu_i^*\in H^1(S^3-L),$$ which we
  can identify with $j^*(\U)$ under $j\colon
  S^3-\nbd{L}\longrightarrow S^3-C_{\p,\q}(L)$, for any choice of
  $\ell$-tuples $\q=(q_1,...q_\ell)$.  Suppose that there is some
  $h_0\in H_1(S^3-L)$ with the property that $\HFLa(L,h)=0$ for all
  $h\in H_1(S^3-L)$ with $h\neq h_0$ and $\langle
  \theta,h\rangle\geq \langle \theta,h_0\rangle$.  Then, we can find
  arbitrarily large $\ell$-tuples $\q=(q_1,...,q_\ell)$ for which the
  following holds.  Letting 
  $$h_1=j_*(h_0)+\OneHalf\sum_{i=1}^{\ell}
  ((p_i-1)\cm (q_i-1)+p_i\cm \sum_{i\neq j}(p_j-1)\cm \Lk(L_i,L_j))\mu_i',$$
  we have that
  $\HFLa(C_{\p,\q}(L),h)=0$ for all $h\in H_1(S^3-C_{\p,\q}(L))$ with
  $h\neq h_1$ and $\langle \U,h\rangle \geq \langle \U,h_1\rangle$.
  Moreover, $\HFLa(C_{\p,\q}(L),h_1)\cong \HFLa(L,h_0)$.\end{prop}

Proceeding as in~\cite{Hedden}, we draw a Heegaard diagram for large
cables of $L$ starting from a Heegaard diagram for $L$. The proof is then
obtained by inspecting the Heegaard diagram. In fact, before giving
the details of the proof, we describe the Heegaard diagram and establish
some of its basic properties, in three lemmas.
  
Recall that for the Heegaard diagram for $L$, each component $L_i$ of
$L$ corresponds to a pair $w_i$ and $t_i$ of basepoints (here, we use
$t_i$ in place of $z_i$, which we reserve for the cable). After
stabilizing the diagram if necessary, we can arrange that the
following conditions hold: 
\begin{itemize} 
\item For $i=1,...,\ell$, $\beta_i$ represents a meridian for the
  corresponding component of $L$, in the sense that $w_i$ and $t_i$
  lie on a curve $\lambda_i$ which meets $\beta_i$ in a single point,
  and is disjoint from all the other $\beta_j$
\item For $i=1,...,\ell$, $\beta_i$ meets some curve $\alpha_i$
  transversally in a single point $x_i$, and is disjoint from all
  the $\alpha_j$ for $j\neq i$.
\end{itemize}
We denote this Heegaard diagram by $(\Sigma,\alphas,\betas,\ws,\ts)$.
        
\begin{figure}[ht]
  \mbox{\vbox{\epsfbox{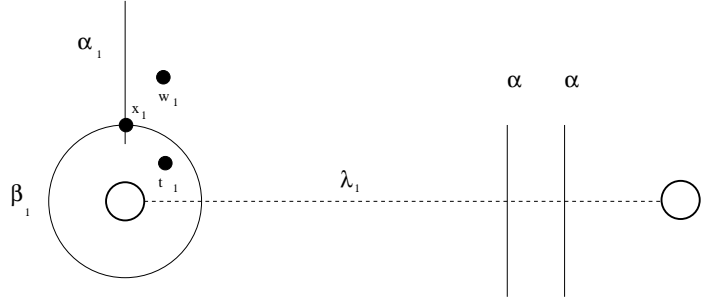}}}
  \caption{\label{fig:PreWinding}
    {\bf{Piece of Heegaard diagram before cabling.}}  After stabilizing a
    Heegaard diagram, we can find a circle $\beta_1$ representing a
    meridian for the first component of a link $L$, so that there is a
    curve $\lambda_1$ which is disjoint from all $\beta_i$ with $i\neq 1$,
    meeting $\beta_1$ in one point.  Note, however, that $\lambda_1$
    typically crosses other $\alpha$-circles, which are indicated here by
    several arcs. The two hollow circles represent a handle to be added to
    the plane. For a more general link, we can find $\ell$ different
    configurations as above.}
\end{figure}

\begin{figure}[ht]
  \mbox{\vbox{\epsfbox{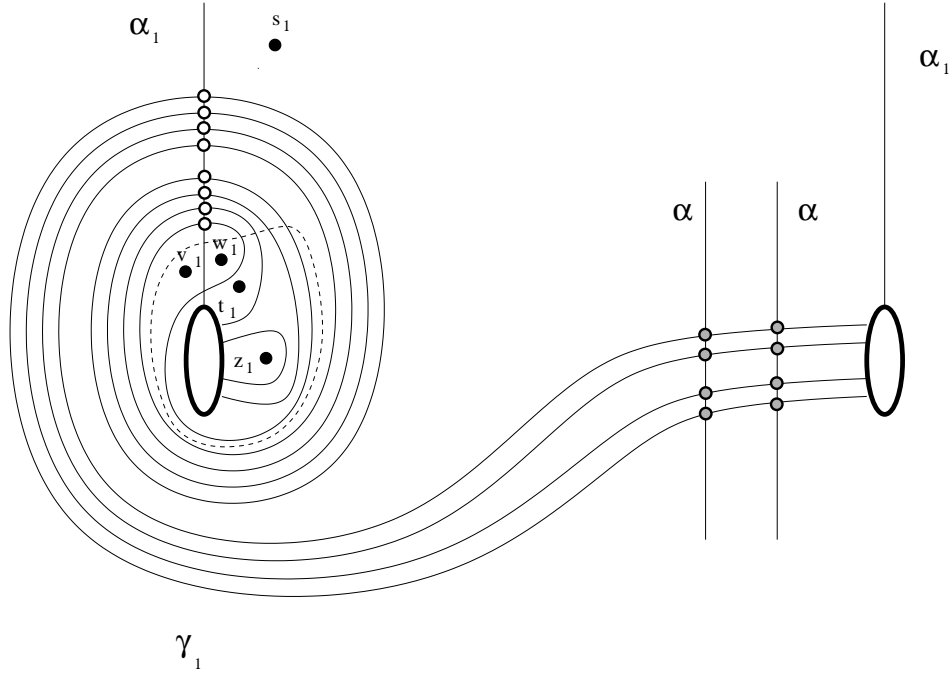}}}
  \caption{\label{fig:Winding} {\bf{Piece of Heegaard diagram after
        cabling.}}  Replace $\beta_1$ from Figure~\ref{fig:PreWinding}
    by a curve $\gamma_1$ which is supported in a neighborhood of
    $\beta_1\cup\lambda_1$.  Possible $\gamma_1$-components of
    exterior intersection points are labelled by the eight hollow
    circles, while possible $\gamma_1$-components of interior
    intersection points are labelled by the eight gray circles. Here,
    $w_1$ and $z_1$ represent the components of the $(3,7)$ cable of
    the component considered in Figure~\ref{fig:PreWinding}. Note that
    if we use the reference point $t_1$ in place of $z_1$, we obtain
    a pointed Heegaard diagram for $L$. The basepoint $s_1$ will be
    used in the proof of Lemma~\ref{lemma:CablingLemma} below. The original
        curve $\beta_i$ (indicated by the dotted line) and the basepoint $v_i$
        will be used in the proof of Lemma~\ref{lemma:OutermostsOnly}.}
\end{figure}

Now, we replace $\beta_i$ with a new curve $\gamma_i$, gotten by
performing a ``finger move'' of $\beta_i$ along $\lambda_i$ with
multiplicity $(p-1)$, and then, in the end, winding some number $n_i$
of times parallel to $\beta_i$.  We then place a new basepoint $z_i$
inside the end of the finger. For notational consistency, we also let
$\gamma_i$ for $i>\ell$ denote the corresponding $\beta_i$. The
resulting diagram $(\Sigma,\alphas,\gammas,\ws,\zs)$ represents the
cable $C_{\p,\q}(L)$, where here
\begin{equation}
  \label{eq:FormOfQ}
  q_i=p_i n_i+1
\end{equation}
for some $\ell$-tuple of integers $\n=(n_1,...,n_\ell)$. 
Note that $(\Sigma,\alphas,\gammas,\ws,\ts)$ still represents
$L$ (the $\gamma_i$ are isotopic to the $\beta_i$ through an isotopy
which crosses only the $z_i$, and no other basepoints).
(Note that we stick with $q_i$ as
in Equation~\eqref{eq:FormOfQ} for concreteness; it is easy to find a
similar description of $C_{\p,\q}(L)$ for other types of $\q$, as
well.)
See Figures~\ref{fig:PreWinding} and \ref{fig:Winding} for an
illustration.

\begin{defn}
Note that $\gamma_i$ is supported in a small regular $N_i$ neighborhood of
$\beta_i\cup \lambda_i$. The intersection points of $\Gen(C_{\p,\q}(L))=
\Ta\cap\Tc$ whose
$\gamma_i$ component is supported in the regular neighborhood of
$\beta_i$ for each $i$ are called {\em $i$-exterior intersection points}
(using terminology of Hedden) and the remaining ones are called  {\em
$i$-interior intersection points}. An intersection point which is 
$i$-exterior for all $i=1,...,\ell$ is called simply an {\em exterior
intersection point}. See Figure~\ref{fig:Winding}.
\end{defn}

In fact, the direction of the winding distinguishes an intersection
point $x_i^0\in\alpha_i\cap\gamma_i$ in the region adjacent to $w_i$.
More specifically, if we consider undoing the finger move (allowing
$\gamma_i$ to cross $z_i$ but not $w_i$ or $t_i$), then the
intersection point of $x_i^0$ corresponds to the original intersection
point $x_i$ between $\alpha_i$ and $\beta_i$. 

\begin{defn}
\label{def:Maximal}
An intersection point $\x'\in\Gen(C_{\p,\q}(L))=\Gen(C_{\p,\q}(L))$
whose $\gamma_i$ coordinate is $x^{0}_i$ for all $i=1,...,\ell$ is
called a {\em maximal exterior} point. Given any $\x\in\Gen(L)$, there
is a corresponding maximal exterior generator
$\x'\in\Gen(C_{\p,\q}(L))$.
\end{defn}

It will be useful to us to
calculate the absolute $\HH$-grading of maximal exterior points.

\begin{lemma}
\label{lemma:HGradingMax}
Fix an intersection point $\x\in \Ta\cap\Tb=\Gen(L)$, and let
$\x'\in\Gen(C_{\p,\q}(L))=\Ta\cap\Tc$ be its corresponding maximal
exterior intersection point whose $\gamma_i$ coordinate (for
$i=1,...,\ell$ is $x_i^0$ (whereas $\x$ has $x_i$ for its $\beta_i$
coordinate). Then,
$$\HGrading_{\w,\ts}(\x)=\HGrading_{\w,\ts}(\x')$$
and also
$$\HGrading_{\w,\z}(\x')=
j_*(\HGrading_{\w,\ts}(\x))+\OneHalf\sum_{i=1}^{\ell}
((p_i-1)\cm (q_i-1)+p_i\cm \sum_{i\neq j}(p_j-1)\cm \Lk(L_i,L_j))\mu_i'.$$
\end{lemma}
  
\begin{proof}
The first claim is easy: the two pointed Heegaard diagrams for $L$,
$(\Sigma,\alphas,\betas,\ws,\ts)$ and
$(\Sigma,\alphas,\gammas,\ws,\ts)$ are isotopic via an isotopy which
does not cross any of the basepoints, and which carries $\x$ to $\x'$.

The second involves  more work. 

First, suppose that $\x,\y\in\Gen(L)$, and $\phi\in\pi_2(\x,\y)$.
Then, it is easy to find $\phi'\in\pi_2(\x',\y')$ which is gotten by
applying a $p_i$-fold finger move to $\phi$ along each of the
$\lambda_i$. For this new domain, we have that
$$n_{z_i}(\phi)-n_{w_i}(\phi)=p_i\cm \left(n_{t_i}(\phi)-
  n_{w_i}(\phi)\right).$$
It follows at once that there is a function
$f(L,\p,\q)$ (independent of $\x\in\Gen(L)$, but depending on the
link $L$; in fact it depends {\em a priori} on the Heegaard diagram we
are using for $L$) with
$$\HGrading_{\w,\z}(\x')-j_*(\HGrading_{\w,\ts}(\x))=
\sum_{i=1}^\ell f_i(L,\p,\q)\mu_i'.$$

Next, we wish to show that for each $i=1,...,\ell$,
$$f_i(L,\p,\q)-\frac{(p_i-1)\cm \sum_{j\neq i}p_j\cm\Lk(L_i,L_j)}{2}$$
is independent of $p_j$ and $q_j$ for $j\neq i$; i.e. there is
a function $\phi_i(p_i,L)$ with the property that
\begin{equation}
\label{eq:DependsOnPj}
f_i(L,\p,\q)=\phi_i(p_i,q_i,L)+\frac{(p_i-1)\cm \sum_{j\neq i}p_j\cm\Lk(L_i,L_j)}{2}.
\end{equation}

The function $f_i(L,\p,\q)$ is understood as follows.  Let $F_i$ be a
Seifert surface for the component $L_i\subset L$, punctured so that it
is supported inside $S^3-\nbd{L}$.  Similarly, let $F_i'$ be a Seifert
surface for the cable $L_i'=C_{p_i,q_i}(L_i)\subset C_{\p,\q}(L)$,
punctured so that it is supported inside $S^3-\nbd{C_{\p,\q}(L)}$.  It
is easy to see from Equation~\eqref{eq:RelSpinCFormula} that
\begin{equation}
\label{eq:InterpretF}
2f_i(L,\p,\q)=\langle c_1(\relspinc_{\w,\z}(\x')),[F_i']\rangle -
p_i \langle c_1(\relspinc_{\w,\ts}(\x),[F_i]\rangle + (p_i-1).
\end{equation}

The intuition behind Equation~\eqref{eq:DependsOnPj} now is the
following. The vector fields determined by $\relspinc_{\w,\z}(\x')$
and $\relspinc_{\w,\ts}(\x')$ agree in $S^3-\nbd{L}$, thought of as a
neighborhood of $S^3-\nbd{C_{\p,\q}(L)}$. Moreover, one can find a
Seifert surface for $L$ in $S^3-\nbd{C_{\p,\q}(L)}$ which has the form
$p_i\cm F_i$ in $S^3-\nbd{L}$. (Note that we are being a bit free with
the meaning of the term {\em Seifert surface}: for our present
purposes, we mean a relative two-dimensional homology class
in the link complement which has intersection number equal to one
with the meridian of $L_i$, and zero with the meridians of all $L_j$ with
$j\neq i$.) Thus,
the difference between the first Chern classes of the two vector
fields over $F_i'$ and $p_i
\cm F_i$ localize to a sum of terms supported near the various
$L_j$. The localization near $L_i$ is independent of $p_j$ for $j\neq
i$, while the local contribution near $L_i$ is $p_i \cm (p_j-1)\cm \Lk
(L_i,L_j)$: this follows from the fact that $F_i'$ meets the $j^{th}$
component of $C_{\p,\q}(L)$ with multiplicity $p_i\cm p_j\cm \Lk
(L_i,L_j)$, whereas $p_i\cm F_i$ meets the $j^{th}$ component of $L$
with multiplicity $p_i\cm \Lk(L_i,L_j)$.

To formulate this intuition precisely, we reformulate the quantities
in terms of the Heegaard diagram. To this end, it is useful to have a
Seifert surface for $L_i\subset L$ drawn on the Heegaard diagram, as
follows.  Let $\xi_i\subset \Sigma-\alpha_1-...-\alpha_{g+\ell-1}$ be
a path from $t_i$ to $w_i$, and $\eta_i\subset
\Sigma-\gamma_1-...-\gamma_{g+\ell-1}$ be another path from $t_i$ to
$w_i$. The closed curve
$\xi_i-\eta_i$ is homologous in $\Sigma$ to a linear combination of
curves chosen among the $\alpha_j$ and the $\gamma_k$ with
$j,k=1,..,g+\ell-1$ but $k\neq i$. Thus, we can find some two chain
$P_i$ in $\Sigma$ representing this homological relation. We assume
without loss of generality (by adding on multiples of regions in
$\Sigma-\alpha_1-...-\alpha_{g+\ell-1}$ if needed) that $P_i$
satisfies $n_{w_j}(P_i)=0$ for $j=1,...,\ell$. 
First, remove disks around the $w_j$ and $t_k$. Next, attach disks to
$P_i$ along the $\alpha_j$ and $\gamma_k$ boundaries.  Finally, attach
a pair of half-disks along the $\xi_i$ and $\eta_i$ arcs. In this
manner, we obtain a Seifert surface $F_i$ for the component
$L_i\subset L$, punctured so as to be supported in the link
complement. 

Note that we can draw $\xi_i$ and $\eta_i$ in the neighborhood $N_i$.
Similarly, we let $\xi_i'$ and $\eta_i'$ be the corresponding paths
with $z_i$ replacing the role of $t_i$. We can construct a two-chain $P_i'$
connecting $\xi_i'-\eta_i'$ in $\Sigma$ with a linear combination of
$\alpha_j$ and $\gamma_k$. We can also build an analogous surface
$F_i'$ for the corresponding component of $C_{\p,\q}(L)$ is obtained
similarly from $P_i'$ by deleting disks around $z_k$.  Clearly, the
two-chains $P_i'$ and $p_i\cm P_i$ are identical, away from the
winding region $N_i$. In particular, both have the same behaviour near
$N_j$ with $j\neq i$, and hence the difference
\begin{equation}
        \label{eq:PartialSeifert}
p_i\cm \langle c_1(\spincrel_{\w,\ts}(\x)),[F_i]\rangle - 
\langle c_1(\spincrel_{\w,\ts}(\x)),[F_i']\rangle  
\end{equation}
is independent of the $p_j$ for $j\neq i$.

There are also Seifert surfaces
$F_i''$ for $L_i$ inside the link $C_{\p',\q}$, where here 
$$p_j'=\left\{\begin{array}{ll}
p_i & {\text{for $i=j$}} \\
1 & {\text{for $i\neq j$}.}
\end{array}
\right.$$
We can draw this on the same Heegaard surface, as follows.
Let ${\mathbf u}$ be the $\ell$-tuple of
points
$$u_j=\left\{\begin{array}{ll}
    z_i' & {\text{if $i=j$}} \\
    t_j & {\text{if $i\neq j$.}}
\end{array}
\right.
$$ The Seifert surface $F_i''$ is obtained from $P_i'$ by puncturing it
in the $t_j$ rather than the $z_j'$.

In fact, it is easy to see that for any $j\neq i$,
$n_{z_j}(P_i)=\lk(L_i,L_j)$. 
Moreover, for fixed intersection point $\x\in\Ta\cap\Tc$,
$\relspinc_{\w,\z}(\x)$ and $\relspinc_{\w,{\mathbf u}}(\x)$ are
represented by the same vector field in $S^3$.  In fact, both
\begin{eqnarray*}
\langle c_1(\relspinc_{\w,\z}(\x)),[F_i']\rangle &{\text{and}}&
\langle c_1(\relspinc_{\w,{\mathbf u}}(\x)),[F_i'']\rangle
\end{eqnarray*}
are obtained by evaluating a relative cohomology class over the
two-chain $P_i'$, appropriately punctured. The difference between
these evaluations comes from the fact that $F_i'$ is obtained by
removing disks $D_j$ around $z'_j$ inside $P_i'$, where the chain
$F_i'$ has multiplicity $\sum_{j\neq i} p_i \cm p_j\cm \lk(L_i,L_j)$,
whereas $F_i''$ is obtained by removing disks around the $t_j$, where
the chain $F_i''$ has multiplicity $p_i\cm \sum_{j\neq i}
\lk(L_i,L_j)$.  Moreover, away from these disks, the two line bundles
associated to $\relspinc_{\w,\z}(\x)$ and $\relspinc_{\w,{\mathbf
    u}}(\x)$ are identified, coming with a canonical trivialization
along $\partial D_j$; whereas along $D_j$, one vector field is gotten
by modifying the other in a prescribed manner (so as to cancel zeros
of $\grad f$, as explained in Subsection~\ref{subsec:RelSpinC}). Hence
the difference is given by
$$(p_i-1)\cm \sum_{j\neq i} p_j\cm \Lk(L_i,L_j) \cm C,$$
where $C$ which depends on the difference between the two
trivializations of the two line fields which extent over the disk. One
can verify that $C=1$ by calculating a model example (the minimal one
being $(2,1)$ cable of the Hopf link).

It follows that 
$$\langle c_1(\relspinc_{\w,\z'}(\x)),[F_i']\rangle
-\langle c_1(\relspinc_{\w,u}(\x)),[F_i']\rangle
= (p_i-1)\sum_{j\neq i} p_j\cm \Lk(L_i,L_j).$$
Combining this with the fact that
$p_i\cm \langle c_1(\spincrel_{\w,\ts}(\x)),[F_i]\rangle - 
\langle c_1(\spincrel_{\w,\ts}(\x)),[F_i']\rangle$ 
is independent of $p_j$ with $j\neq i$
(cf. Equation~\eqref{eq:PartialSeifert}), together with 
the interpretation of $f_i$ from Equation~\eqref{eq:InterpretF},
Equation~\eqref{eq:DependsOnPj} follows.

Next, we consider the 
dependence of $\phi_i(p_i,q_i,L)$ on $q_i$.  If we fix
$p_i$, then the quantity $$\phi_i(p_i,q_i+p_i,L)-\phi_i(p_i,q_i,L)$$
localizes around $N_i$, and is independent of $L$.
This is true since, once again, the two-chains representing the Seifert
surfaces and the vector fields representing corresponding intersection points
differ  only near $N_i$.

By considering a model calculation, one can see that
$$\phi_i(p_i,q_i,L)=\frac{q_i (p_i-1)}{2}+\psi_i(L,p_i).$$ The
simplest model calculation, of course, is a $(p_i,p_i n_i+1)$-cable of
the unknot (endowed with a genus one Heegaard diagram with a single
generator $x$). In this case, the Heegaard diagram described above is
a diagram for the $(p_i,q_i)$ torus knot, with at most one generator
in each $\HH$-grading. It is straightforward to see that $x_0$, here
is the generator with maximal $\HH$-grading, which, of course, then
agrees with the highest $T$-power of the (symmetrized) Alexander
polynomial, $\frac{(p_i-1)(q_i-1)}{2}$.

In a similar manner, if we vary $p_i$, we have that
$$\phi_i(p_i,q_i,L)=\frac{(p_i-1)(q_i-1)}{2}+c(L).$$
Obviously, setting $p_i=1$, we see that $c(L)=0$.
\end{proof}

For $i=1,...,\ell$,
$\alpha_i\cap\gamma_i$ consists of $2 n_i (p_i-1)+1$ intersection
points. 
Given $x,x'\in \alpha_i\cap\gamma_i$, it is easy to see that there
are arcs $a\subset \alpha_i$ and $b\subset \gamma_i$, both going
from $x$ to $x'$, with the additional property that $a-b$ is homologous
to a sum of curves among the $\alpha_m$ and $\beta_n$. Let $D_{x,x'}$
be such a two-chain.
Consider the function
$$\zeta_j\colon \alpha_i\cap\gamma_i\longrightarrow \Z$$ which
is uniquely characterized up to overall translation by the equation
$$\zeta_j(x)-\zeta_j(x')=n_{z_j}(D_{x,x'})-n_{w_j}(D_{x,x'}).$$
As an
immediate consequence of Equation~\eqref{eq:DefHGrading}, we see that
if $\x,\x'\in \Ta\cap\Tc$ are two intersection points which agree on
all factors except for $i=1,...,\ell$ on
$\alpha_i\cap\gamma_i$, where $\x$ is $x_i$ and
$\x'$ is $x'_i$, then
$$\HGrading_{\w,\ts}(\x)-\HGrading_{\w,\z}(\x')=\sum_{i=1}^{\ell} (\zeta_i(x_i)-\zeta_i(x_i'))\cm \mu_i.$$

\begin{figure}[ht]
  \mbox{\vbox{\epsfbox{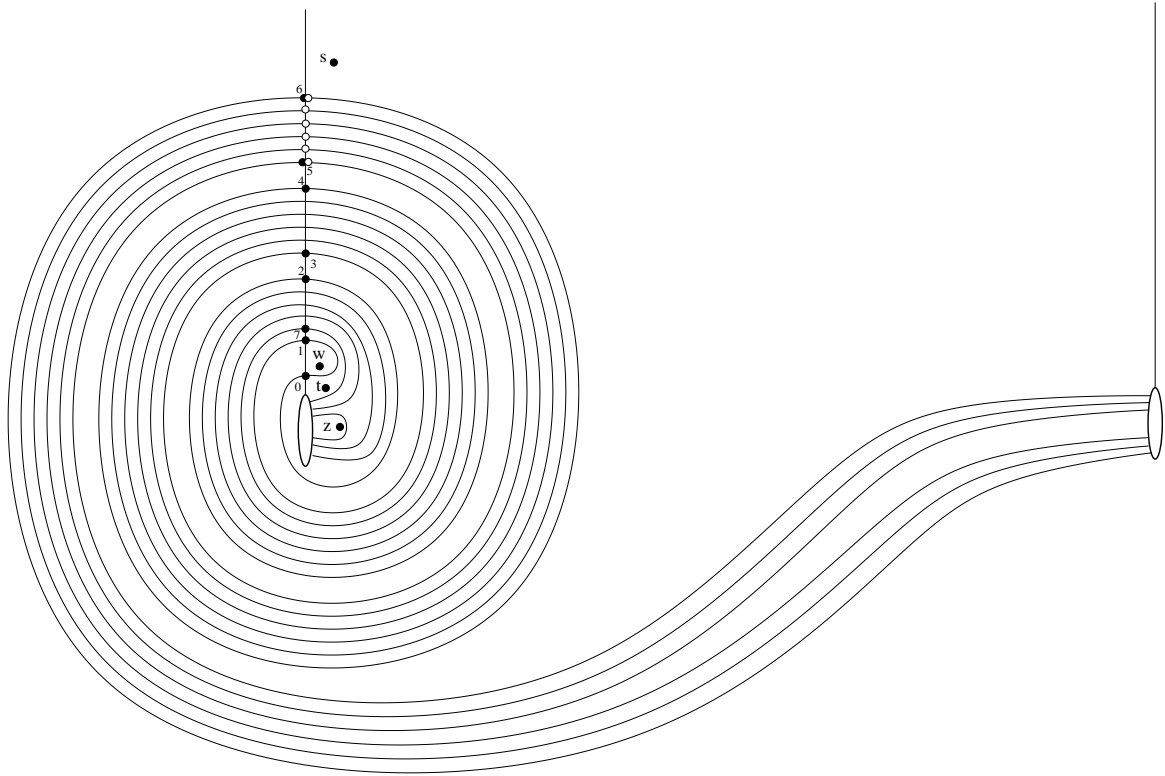}}}
  \caption{\label{fig:Cable} {\bf{Cable with $p=4, q=13$.}}
    We have dropped all subscripts. We have indicated the $x^i_1$ for
    $i=0,...,7$ (but listed them only by $i$); those with $i\leq 6$
    are outermost. Innermost exterior points are also indicated (but
    not labelled) with hollow circles (though note that $i=5$ and $6$
    are both innermost and outermost.}
\end{figure}


\begin{lemma}
\label{lemma:ZetaProperties}
We can order the intersection points of $\alpha_i\cap\gamma_i$
$\{x_i^k\}_{k=0}^{2 (p_i-1) n_i}$ with the convention that
$\zeta_i(x_i^j)>\zeta_i(x_i^k)$ if $j<k$.
\begin{eqnarray}
  \label{eq:Ordering}
  \zeta_i(x_i^j)>\zeta_i(x_i^k)&{\text{if}}& j<k.
\end{eqnarray}
For the function $\zeta_i$ as above, we have that
\begin{equation}
  \label{eq:Breadth}
  \zeta_i(x^{0}_{i})-\zeta_i(x^{2n_i}_{i})=p_i n_i.
\end{equation}
Moreover, for
$i\neq j$  we have that
\begin{equation}
  \label{eq:SmallJunk}
\zeta_j(x^{k}_i)-\zeta_j(x^{k+1}_{i})
=\left\{
\begin{array}{ll}
p_j\cm \lk(L_i,L_j)& \text{if $2n_i|k$} \\
0 &{\text{otherwise.}}
\end{array}
\right.
\end{equation}
\end{lemma}

\begin{proof}
  Equation~\eqref{eq:Breadth} can be verified by constructing domains
  which are supported entirely inside $N_i$.  Starting from 
  $x^{0}_i$ as in Definition~\ref{def:Maximal}, we can
  define $x^m_i$ for $m=0,...,2n-2$,
  in such a manner that there is a bigon from 
  $x_i^{2k}$ to $x_i^{2k+1}$, supported
  in $N_i$, with local multiplicity $-1$ at $w_i$ (and multiplicity
  zero at $z_i$), and an immersed bigon connecting $x_i^{2k+1}$ to
  $x_i^{2k+2}$, supported in $N_i$ with local multiplicity $(p-1)$ at
  $z_i$ (and multiplicity zero at $w_i$), where here $0\leq 2k \leq 2p-2$.
  (See Figure~\ref{fig:DomainsInCable} for an illustration.)
  Adding these up, we get Equations~\eqref{eq:Breadth}.

  \begin{figure}[ht]
    \mbox{\vbox{\epsfbox{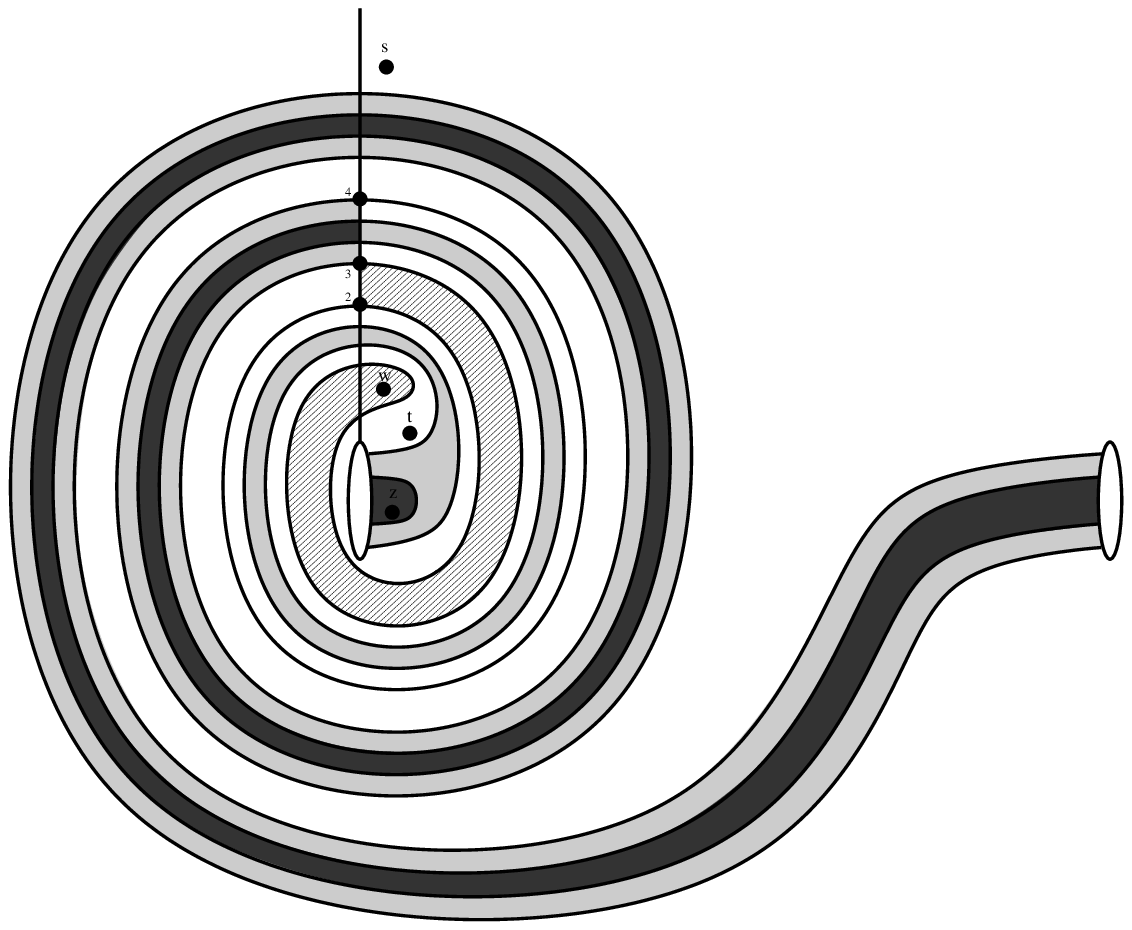}}}
    \caption{\label{fig:DomainsInCable}
      {\bf{Domains illustrating Equation~\eqref{eq:Breadth}.}}  There is
      a domain from $x^2$ to $x^3$ (indicated here by the hatched line)
      which crosses $w$ with multiplicity $-1$, and one from $x^3$ to
      $x^4$ which crosses $z$ with multiplicity $2$.  (Regions with
      multiplicity $+1$ are shaded light gray, those with $+2$ are
      shaded dark gray.)}
  \end{figure}

  Indeed, in a similar manner, we can extend the ordering so that
  there is domain (always an immersed disk) connecting $x_i^k$ and
  $x_i^{k+1}$ which is supported entirely inside $N_i$, provided that
  $2n_i$ does not divide $k$. In particular, it follows that when
  $2n_i$ does not divide $k$,
  $\zeta_j(x_i^k)-\zeta_j(x_i^{k+1})=0$. Moreover, with these
  conventions (and depending on the parity of $k$) the disk always has
  either multiplicity $0$ at $w_i$ and positive at $z_i$, or it has
  multiplicity $0$ in $z_i$ and negative multiplicity at $w_i$.

  In the special cases where $2n_i|k$, however, there is a domain
  whose boundary contains the part of $\alpha_1$ outside the spiral
  region. In completing this domain, we find that
  $\zeta_j(x_i^k)-\zeta_j(x_i^{k+1})$ is given by $p_j\cm
  \lk(L_i,L_k)$.
\end{proof}

Note that that $x^0_i$ is the $\gamma_i$-coordinate of a maximal exterior
intersection point in the sense of Definition~\ref{def:Maximal}.

\begin{defn}
  \label{def:Exterior}
  An exterior intersection point $\x$ whose $\gamma_i$ coordinate 
        satisfies $0\leq k\leq 2(n_i-1)$ according to the above
        labeling convention, then we say that $\x$ is $i$-outermost.
  An exterior intersection point $\x$ whose $\gamma_i$
        cordinate is one of the $2(p_i-1)$ points among the $x_i^k$
        closest to $s_i$
        (i.e. one of $x_i^k$ where $k$ is of $2(p_i-1)$ possible integers with 
        $1\leq k\leq 2(p_i-1)n_i$ with $2n_i|k$ or $2n_i|k+1$
        is called {\em $i$-innermost}.
\end{defn}

Note that an exterior intersection point can be both $i$-innermost and
$i$-outermost at the same time (if its $\gamma_i$ component is
either $2n_i-1$ or $2n_i$).

Recall that a maximal exterior intersection point in the sense of
Definition~\ref{def:Maximal} is an outermost intersection point.  See
Figure~\ref{fig:Cable} for an illustration.

Sometimes, when stressing the dependence of $\zeta_i$ on the winding parameter
$\n$, we write $\zeta_i^{\n}$.

\begin{lemma}
  \label{lemma:OutermostsOnly} Given $\p=(p_1,...,p_\ell)$ and $\ell$-tuple
  of positive integers, each of which is greater than one. Then, all
  sufficiently large $\n=(n_1,...,n_\ell)$ have the following
  property.  For $\q=(q_1,...,q_\ell)$ as in
  Equation~\eqref{eq:FormOfQ}, if $\x\in\Gen(C_{\p,\q}(L))$ is a
  generator with the property that $$\langle \U,
  \HGrading_{\w,\z}(\x)\rangle \geq \langle \U,
  \HGrading_{\w,\z}(\x_0')\rangle $$
  for some $\x_0\in\Gen(L)$, then
  $\x$ is an outermost exterior point.
\end{lemma}  

\begin{proof}
We claim that there is a constant $c=c(\p,L)$ with the property that
for all sufficiently large $\n$, if $\y$ is an intersection point
which is not $i$-exterior, then there is an maximal exterior point
$\x=\x(\n)$ with the property that
\begin{equation}
  \label{eq:OneConstant}
\langle\U,\HGrading_{\w,\z}(\x)\rangle -
\langle\U,\HGrading_{\w,\z}(\y)\rangle
\geq c+ n_i.
\end{equation}
The lemma will follow at once from this, together with the fact
that there is a universal bound $c'$ independent of $\n$ with the property that if
$\x,\y\in\Gen(L)$, then 
$$|h_{\w,\z}(\x')-h_{\w,\z}(\y') |<c'.$$
(This latter fact follows at once from Lemma~\ref{lemma:HGradingMax}.)

Equation~\eqref{eq:OneConstant} is established as follows. Suppose
that $\n_0$ is fixed, and fix some intersection point $\y_0$ which is
$i$-interior for some $i$.  We claim indeed that there is an exterior
intersection point $\x_0$ which we connect to $\y_0$ via a domain
$\phi\in\pi_2(\x_0,\y_0)$ whose multiplicity at $w_i$ is zero. In
fact, we claim that for a suitable such choice, we can arrange also
that the local multiplicity of $\phi$ also at $v_i$ is zero (where here
$v_i$ is the basepoint separated from $w_i$ by an arc crossing only
$\alpha_i$, as pictured in Figure~\ref{fig:Winding}).  We argue this
as follows. Choosing $\x_0$ to be $i$-innermost, we can arrange that its
$\gamma_i$ coordinate $x_i^k$ can be connected to the
$\gamma_i$-coordinate of $\y_0$ by an arc which meets no other $x_i^j$
for $j\neq k$. Similarly, we can connect $x_i^k$ to the $\alpha_1$
coordinate of $\x_0$ by an arc which points out of the winding region
(an hence meeting only innermost intersection points with
$\gamma_i$). In particular, both arcs are disjoint from the dotted
curve ($\beta_i$) indicated in Figure~\ref{fig:Winding}.  We complete
these two arcs to a choice of curves $\epsilon$ composed of arcs among
the $\alpha_i$ and $\gamma_j$, where the transitions alternate between
points in $\x_0$ and points in $\y_0$. Any domain connecting $\x_0$
and $\y_0$ is gotten from $\epsilon$ by adding sufficient multiples of
the $\alpha_k$ and $\gamma_\ell$ to make it null-homologous. We claim
that in this procedure, there will be no copies of $\alpha_1$
added. This is clear since the algebraic intersection number of the
original dotted curve $\beta_i$ with all other curves is zero, and also
our original curve $\epsilon$ does not cross $\beta_i$. It follows now
that the multiplicity of our domain at $w_i$ agrees with its
multiplicity at $v_i$.

Increasing the winding parameter to $\n$, we we claim that there is
always some exterior intersection point $\x_\n$ with the following properties.
The $\beta_k$ with $k>\ell$ components of $\y_n$ are fixed,
coinciding with those for $\y_{\bf 0}$;
for each $i$ with the property that $\y_{\bf 0}$ is $i$-interior,
$\x_\n$ is $i$-innermost, in the sense of Definition~\ref{def:Exterior}.
Also, there is a homotopy class
$\phi_\n\in\pi_2(\x_\n,\y_\n)$ with the property that
\begin{eqnarray*}
\zeta_i^{\n}(\x_\n)-\zeta_i^{\n}(\y_\n)
&=&
n_{w_i}(\phi_\n)-n_{z_i^{\n}}(\phi_\n)\\
&=&n_{w_i}(\phi)-n_{z_i^{\n_0}}(\phi) \\
&=&\zeta_i^{\n_0}(\x_0)-\zeta_i^{\n_0}(\y).
\end{eqnarray*}
(Note we are using here basepoints $z_i^{\n}$ and $z_i^{\n_0}$ for two
different $\ell$-tuples of winding parameters $\n$ and $\n_0$; and we
record this in the notation for $\zeta_i$.)  To see this, note that
the new curves $\gamma_i^{\n}$ (for $i=1,...,\ell$) are gotten by
performing Dehn twists to $\gamma_i^{\n_0}$ along the curves $\beta_i$
(again, as in Figure~\ref{fig:Winding}).  The new domain $\phi_\n$ is
gotten from the original domain $\phi$ by performing Dehn twists along
its $\gamma_i$-boundary. (The domain can then be used to determine the
intersection point point $\x_\n$.). This Dehn twist can be done to
obtain a new domain precisely since the multiplicities at $v_i$ and
$w_i$ agree.

Now, for each $j$ for which  $\y_{\bf 0}$ is $j$-exterior, 
since $\phi^\n$ is gotten from $\phi$ by a local procedure near $N_i$,
we still have have that
\begin{eqnarray*}
\zeta_j^{\n}(\x_\n)-\zeta_j^{\n}(\y_\n)
&=&n_{w_j}(\phi_\n)-n_{z_j^{\n}}(\phi_\n) \\
&=&n_{w_j}(\phi)-n_{z_j^{\n}}(\phi) \\
&=&\zeta_j^{\n_0}(\x_0)-\zeta_j^{\n_0}(\y)
\end{eqnarray*}
is bounded independent
of $\n$.
Combining 
Equations~\eqref{eq:Ordering}, \eqref{eq:Breadth}, and~\eqref{eq:SmallJunk}
we see that
$$\zeta_i(\x'_\n)-\zeta_i(\y_\n)=c_1+ \zeta_i(\x'_n)-\zeta_i(\x_\n)\geq p_i n_i+c_2$$
and for $j$ for which $\y$ is $j$-exterior
$$\zeta_j(\x'_\n)-\zeta_j(\y_\n)\geq c_3+ \zeta_j(\x'_n)-\zeta_j(\x_\n)\geq c_3$$
where here the constants $c_1$, $c_2$, and $c_3$ depend  on only
the interior
part of $\y_0$ (and are in particular independent of $n_i$).  Since
there is only a finite number of possibilities for this interior part
of any intersection point $\y$, we can find one constant 
as required in Inequality~\eqref{eq:OneConstant}.
\end{proof}

Consider the map $\HH(L)\longrightarrow \HH(C_{\p,\q}(L))$
which carries
$\HGrading_{\w,\ts}(\x)\longrightarrow \HGrading_{\w,\z}(\x')$,
where here $\x'$ is the maximal exterior point nearest to $\x$.
We denote the induced map by $h\mapsto h'$.

Given $h,k\in \HH(L)$, we write $k\geq h$ if 
$$k=h+\sum_{i=1}^\ell a_i\cm \mu_i,$$ where $a_i$ are non-negative integers.

\begin{lemma}
  \label{lemma:CablingLemma} Fix a link $L$ and a cabling parameters
  $\p$ and $\q$ as in Lemma~\ref{lemma:OutermostsOnly}. Fix $h\in
  \HH(L)$ and suppose moreover that $\HFLa(L,k)=0$ for all $k\geq
  h$. Then, if $h'$ is represented by outermost exterior intersection
  points only, we have that $\HFLa(C_{\p,\q}(L),h')\cong \HFLa(L,h)$.
  Moreover, if $k\in \HH(C_{\p,\q}(L))$ is represented entirely of
  outermost intersection points and $\HFLa(C_{\p,\q}(L),k)\neq 0$,
  then there is some $h\in\HH(L)$ with $h'\geq k$ and $\HFLa(L,h)\neq
  0$.
\end{lemma}

\begin{proof}
  This will follow from a spectral sequence whose $E_2$ term consists
  of
  $$\bigoplus_{k\in S} \HFLa(L,k),$$
  converging to $\HFLa(C_{\p,\q}(L),h')$, where here $S\subset \HH(L)$
  is some set of $k$ with $k\geq h$, and which contains $h$.

  The spectral sequence is constructed using additional basepoints
  $\sss=\{s_1,...,s_\ell\}$ placed outside the winding region, as in
  Figure~\ref{fig:Winding}. These basepoints induce an additional
  filtration on $C_{\p,\q}(L,h')$. We can think of this filtration
  concretely in the following terms: for exterior intersection points,
  the generators in a fixed $s_i$-filtration are those whose $i^{th}$
  component is some fixed intersection point $x^i_j$. The homology of
  the associated graded object counts holomorphic disks which do not
  cross the spiral region. For fixed $h'$, we can make the filtration
  $\Z$-valued (rather than relative) by the convention that maximal
  intersection points $\x'$ have $s_i$-filtration equal to zero.
  Indeed, we will put these filtrations together, and consider the
  filtration by $\m=(m_1,...,m_\ell)$. Let $\sigma^{\m}
  \CFLa(C_{\p,\q}(L),h')$ be the associated graded complex, generated
  by $\x$ whose $s_i$ filtration is given by $m_i$.
  
  We claim that 
  \begin{equation}
    \label{eq:CableIsomorphism}
    H_*(\sigma^{\mathbf{0}} \CFLa(C_{\p,\q}(L),h'))\cong
    \HFLa(L,h).
  \end{equation}
  For this, we consider the $\sss$ as inducing a
  filtration on the chain complex $\CFLa(L,h)$, gotten by using the
  basepoints $\ws$ and $\ts$. Observe that for this complex, we can
  isotope the $\gamma_i$ to $\beta_i$ (crossing the $\zs$ but none of
  the other basepoints), so that in the end the basepoints $s_i$ is in
  the same component as $w_i$. Thus, it induces a trivial filtration.
  We consider the $\zs$ as giving a further filtration, denoted by
  $\zeta$, on $\bigoplus_{\m} \sigma^\m \CFLa(L,h)$.  The homology
  groups are supported entirely inside $\sigma^{\mathbf 0}
  \CFLa(L,h)$: if $\m\neq 0$, then there are bigons preserving
  elements of $\sigma^{\m}\CFLa(L,h)$ with positive multiplicity on
  the $\zs'$ which cancel generators in pairs.  It follows that
  \begin{equation}
    \label{eq:CableIsomorphism0}
    H_*(\sigma^{\mathbf{0}} \CFLa(L,h))\cong \HFLa(L,h).
  \end{equation}
  There is an
  also an easily seen identification of chain complexes
  $$\sigma^{\mathbf{0}} \CFLa(C_{\p,\q}(L),h')\cong
  \sigma^{\mathbf{0}} \CFLa(L,h).$$
  (Here we are using the fact that the equivalence class of $h'$
  uses only outermost generators, for which the differentials then coincide.)
  which, together with Equation~\eqref{eq:CableIsomorphism0},
  gives Equation~\eqref{eq:CableIsomorphism}.
  
  We claim also that for each $\m>0$, $H_*(\sigma^\m C_{\p,\q}(L,h'))$
  is identified with $\HFLa(L,k)$, for some $k>h$. In this case, counting
  differentials crossing $w_i$ once gives an identification 
  \begin{equation}
    \label{eq:CableIsomorphism2}
    H_*(\sigma_i^{2m+1}(C_{\p,\q}(L),h'))\cong 
    H_*(\sigma_i^{2m} (C_{\p,\q}(L),h'+\mu'_i)),
  \end{equation}
  while counting
  differentials crossing $z_i'$ (with multiplicity $p_i$) gives an
  identification 
  \begin{equation}
    \label{eq:CableIsomorphism3}
    H_*(\sigma_i^{2m_i+1}(C_{p,q}(L),h'))\cong 
    H_*(\sigma_i^{2m_i+2}(C_{p,q}(L),h'+p_i\cm \mu'_i))
  \end{equation}
  for all $m_i\geq 0$ (these were the domains used to establish
  Equation~\eqref{eq:Breadth}; again, we are using the fact that
  equivalence classes contain only outermost generators).  It follows now from
  Equations~\eqref{eq:CableIsomorphism}, \eqref{eq:CableIsomorphism2},
  and~\eqref{eq:CableIsomorphism3} together that if $\m>0$, then
  $$H_*(\sigma^{\m}\CFLa(C_{\p,\q}(L),h'))=0.$$ Thus,
  $$H_*(\CFLa(C_{\p,\q}(L),h'))\cong H_*(\sigma^{\mathbf
    0}\CFLa(C_{\p,\q}(L),h'),$$
  and hence applying
  Equation~\eqref{eq:CableIsomorphism}, we obtain the desired
  identification $$\HFLa(C_{\p,\q}(L),h')\cong \HFLa(L,h).$$  Similarly,
  Equations~\eqref{eq:CableIsomorphism}, \eqref{eq:CableIsomorphism2},
  and~\eqref{eq:CableIsomorphism3} together give the second claim.
\end{proof}

\begin{proof} [Of Proposition~\ref{prop:Cabling}.]
  Choose $\n$ large enough as required by
  Lemma~\ref{lemma:OutermostsOnly}.  Given $h_0\in\HH(L)$ and
  $h_1\in\HH(C_{\p,\q}(L))$ as in the statement of the proposition, we
  can also consider $h_0'$, which is represented by maximal
  intersection points corresponding to generators from $h_0$.  Then,
  according to Lemma~\ref{lemma:HGradingMax}, $h_1$ coincides with
  $h_0'$. Moreover, according to Lemma~\ref{lemma:CablingLemma},
  $\HFLa(C_{\p,\q}(L),h_1)\cong \HFLa(L,h_0)$. Indeed, suppose that
  $k\in \HH(C_{\p,\q}(L))$ satisfies $\HFLa(C_{\p,\q}(L),k)\neq 0$ and
  $\langle \U,k\rangle\geq \langle \U,h_1\rangle$. Again, according to
  Lemma~\ref{lemma:CablingLemma}, we have some $h$ with $h'\geq k$ and
  $\HFLa(L,h)\neq 0$. It follows that $\langle \U,h'\rangle \geq
  \langle \U,h_0'\rangle$. We conclude that $h_0=h$, and hence that
  $k=h_0'$, as required.
\end{proof}

Proposition~\ref{prop:Cabling} has the following immediate consequence:

\begin{lemma}
  \label{lemma:FloerNormCable}
  Fix an oriented link $L$ and also an $\ell$-tuple of positive
  cabling coefficients $\p=(p_1,...,p_\ell)$, each of which is greater
  than one.  There are arbitrarily large $\q=(q_1,...,q_\ell)$ with
  $p_i$ and $q_i$ relatively prime, so that the following relation
  holds between the link Floer homology norms of $L$ and
  $C_{\p,\q}(L)$. Consider the homology class $\U\in
  H^1(S^3-C_{\p,\q}(L))$ given by $\U(\mu_i')=1$ for all
  $i=1,...,\ell$. Then,
  $$y(C_{\p,\q}(L),\U)=y(L,j^*(\U))+\sum_{i=1}^\ell \left(\frac{(q_i-1-Q_i)\cm
  (p_i-1)}{2}\right).$$
\end{lemma}

\begin{proof}
This is an immediate consequence of Proposition~\ref{prop:Cabling}.
\end{proof}

\vskip.2cm
\noindent{\bf Proof of Theorem~\ref{thm:ThurstonNorm}.}
It suffices to verify Theorem~\ref{thm:ThurstonNorm} for $h\in H^1(S^3-L;\Z)$ determined by 
$\langle h,\mu_i\rangle = p_i$ for $i-1,...\ell$ and satisfying the following 
two conditions:
\begin{list}
        {(C-\arabic{bean})}{\usecounter{bean}\setlength{\rightmargin}{\leftmargin}}
\item 
\label{n:NowhereVanish}$|p_i|>1$ for all $i=1,...,\ell$
\item
\label{n:GenericMax} there is a unique 
$s\in \HH$ with $\HFLa(L,s)\neq 0$ and maximal evaluation
$\langle s, h\rangle$.
\end{list}
We can see this as follows. Let $M$ denote the set of $h$ with the
above two properties. Clearly, the set of points in $H^1(S^3-L;\R)$
with the property that $r\cm h\in M$ for some $r\in\R$ forms a dense
set. Since both $x(\PD[h]) +\sum_{i=1}^\ell |\langle h, \mu_i\rangle
|$ and $2y(h)$ are continuous functions which are linear on rays, it
suffices to verify that the coincide for elements on $M$.

Without loss of generality, we can orient $L$ so that $p_i>0$ for each
$i=1,...,\ell$.  For $\p=(p_1,...,p_\ell)$,
we can realize $h=j^*(\U)$ for any cable
$C_{\p,\q}(L)$. Indeed, we can make $\q$ arbitrarily large, so that
both Lemmas~\ref{lemma:ThurstonNormCable}
and~\ref{lemma:FloerNormCable} hold.  Then, according to
Lemma~\ref{lemma:FloerNormCable},
\begin{eqnarray*}
2y(L,j^*(\U))&=&2y(C_{\p,\q}(L),\U) - \sum_{i=1}^\ell (q_i-Q_i-1)\cm
(p_i-1).
\end{eqnarray*}
Indeed, by Lemma~\ref{lemma:FloerNormCable} and Condition~\ref{n:GenericMax},
Proposition~\ref{prop:NiPlusSS} applies, to show that 
$$2y(C_{\p,\q}(L),\U)= x(C_{\p,\q}(L),\U)+\ell.$$
Combining this with Lemma~\ref{lemma:ThurstonNormCable}, we see that
\begin{eqnarray*}
2y(L,j^*(\U))&=& x(C_{\p,\q}(L),\U) + \ell - \sum_{i=1}^\ell (q_i-Q_i-1)\cm
(p_i-1) \\
&=& x(L,j^*(\U)) + \sum_{i=1} p_i.
\end{eqnarray*}

We have verified Theorem~\ref{thm:ThurstonNorm} for $h\in
H_1(S^3-L;\Z)$ satisfying Conditions~(C-\ref{n:NowhereVanish}) and
(C-\ref{n:GenericMax}).  Since the set of such homology classes is the
complement of finitely many hyperplanes,
Theorem~\ref{thm:ThurstonNorm} is easily seen to follow for all $h\in
H_1(S^3-L;\R)$.
\qed

\vskip.2cm
\noindent{\bf Proof of Corollary~\ref{cor:AltLinks}.}
According to Theorem~\ref{thm:AltLinks}, the Newton polytope of the
multi-variable Alexander polynomial of an alternating link, when added
to a (suitably centered) unit hypercube, gives a polytope which
can then be scaled by a factor of two to obtain the link Floer homology 
polytope. The result is now an immediate consequence of
Theorem~\ref{thm:ThurstonNorm}.  \qed

\section{On fibered links}
\label{sec:Fibrations}

A direction $\theta\in H^1(S^3-L;\Z)$ is said to be {\em fibered} if it can
be represented by a fibration. Explicitly, this means that there is a
nowhere vanishing one-form $\omega$ defined over $S^3-\nbd{L}$
representing the given cohomology class $\theta$ whose restriction
to $\partial \nbd{L}$ also vanishes nowhere.  An orientation on $L$
gives rise to a canonical cohomology class $\U\in H^1(S^3-L;\Z)$,
whose value on each (oriented) meridian is one.  An oriented link
$\orL$ is fibered if its corresponding cohomology  class
$\U\in H^1(S^3-L;\Z)$ is fibered.  It is easy to see that if $L$ is a
link and $\theta\in H^1(S^3-L;\Z)$ is a fibered cohomology class, with
$p_i=\theta(\mu_i)$, then $C_{\p,\q}(L)$ is a fibered link.

In~\cite{HolDiskContact}, it is shown that if a knot $K\subset Y$ is a
null-homologous knot in a three-manifold which is also fibered, then
its topmost non-trivial Floer homology group is one-dimensional.
Suppose now $\orL$ that is an oriented link, then
it is easy to see that 
$\kappa(L)$ is a fibered knot in $\#^{\ell-1}(S^2\times S^1)$.
Indeed, if $S^3-\nbd{L}$ is the mapping torus of an automorphism
$\phi$ of an oriented surface-with-boundary $F$ with $\ell$ boundary
components, 
and $F'$ denotes the surface obtained by attaching $\ell-1$ one-handles
to $F$ to get a surface with connected boundary, 
then $\#^{\ell-1}(S^2\times S^1)-\kappa(L)$ is the mapping torus 
of the automorphism $F'$ obtained by extending $\phi$ by the
identity map on the new one-handles, see
also~\cite{Ni}.  Thus, from the statement for knots, it follows at
once that the top-most non-trivial knot Floer homology 
$\HFKa(\orL,s)$ of an oriented, 
fibered link $\orL$ is one-dimensional.

\vskip.3cm
\noindent{\bf{Proof of Proposition~\ref{prop:Fibrations}.}}
Let $B_T\subset H^1(S^3-L;\R)$ denote the Thurston polytope and
$B_T^*\subset H_1(S^3-L;\R)$ its dual polytope. Let $Q$ be the symmetric 
hypercube in $H^1(S^3-L)$ with edge-length two.

Fix an extremal point $P$ in $B_T^*$, which corresponds to a fibered face
in $B_T$. This means that there is a cohomology class $\theta\in
H^1(S^3-L;\Z)$ belonging to a fibration of $S^3-L$, with the property
that $\theta(P)\geq \theta(h)$ with equality only if $h=P$.  

It is now easy to find, for each $h_0\in s(P)$, a fibered cohomology
class $\theta_0\in H^1(S^3-L;\Z)$ with the above property (i.e. that
$\theta_0(P)\geq \theta_0(h)$ with equality only if $h=P$) and which
satisfies the additional property that $\theta_0(h_0)\geq \theta(h)$ for
all $h\in B_T^*+Q$ with equality only if $h=h_0$. 
In view of Theorem~\ref{thm:ThurstonNorm}, this additional property 
is equivalent to the
condition that
\begin{equation}
\label{eq:ExtremalCase}
\begin{array}{llll}
\HFLa(L,h)=0 &{\text{if $\theta_0(h)\geq \theta_0(h_0)$ and $h\neq h_0$}}, 
&{\text{whereas}} &
\HFLa(L,h_0)\neq 0.
\end{array}
\end{equation}

Given $h_0$, we can find $\p$ with the property that for any $\q$,
$h_0=j^*(\U)$ for any cable $C_{\p,\q}(L)$. 
Proposition~\ref{prop:Cabling} can be combined
with the above to show that if
we write 
$$h_0'=
j_*(h_0)+\OneHalf\sum_{i=1}^{\ell}
  ((p_i-1)\cm (q_i-1)+p_i\cm \sum_{i\neq j}(p_j-1)\cm \Lk(L_i,L_j))\mu_i',$$
then 
\[\begin{array}{ll}
\HFLa(C_{\p,\q}(L),h')=0 
&{\text{if $\U(h)\geq \U(h_0')$ and $h'\neq h_0'$, whereas}},\\
\HFLa(C_{\p,\q}(L),h_0') \cong \HFLa(L,h_0). &
\end{array}\]
Combining this with
Lemma~\ref{lemma:HFKandHFL}, we conclude that the top-most non-trivial
knot Floer homology group of $C_{\p,\q}(L)$ is isomorphic to
$\HFLa(L,h_0)$. On the other hand, since $h_0$ is a fibered direction,
it follows  that $C_{\p,\q}(L)$
is a fibered link, and hence we conclude that
this homology group $\HFLa(C_{\p,\q}(L),h_0') \cong \HFLa(L,h_0)$ is 
one-dimensional, as claimed.
\qed

\section{Examples}
\label{sec:Examples}

\subsection{An alternating example: $9^2_{41}=9a42$}

\begin{figure}[ht]
\mbox{\vbox{\epsfbox{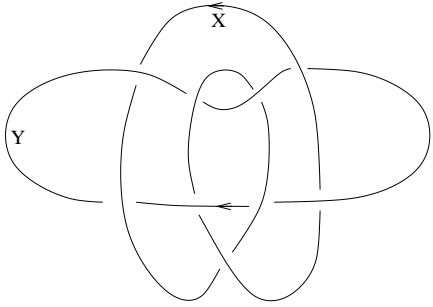}}}
\caption{\label{fig:9a42}
{\bf{The oriented link $9a{42}$.}}
We have specified orientations on the two components,
and labeled them by the symbols $X$ and $Y$.}
\end{figure}

Consider the alternating knot $9^2_{41}$ from Rolfsen's
table~\cite{Rolfsen}, or $9a{42}$ in Thistlethwaite's notation, cf.
Figure~\ref{fig:9a42}.  The symmetrized Alexander of this link is
given by
$$
\begin{array}{rrrr}
-X^{-\frac{3}{2}} Y^{\frac{3}{2}} & + X^{-\frac{1}{2}} Y^{\frac{3}{2}} & & \\
+2  X^{-\frac{3}{2}} Y^{\frac{1}{2}} &  -5  X^{-\frac{1}{2}} Y^{\frac{1}{2}} &     +4 X^{{\frac{1}{2}}} Y^{\frac{1}{2}}&  - X^{{\frac{3}{2}}} Y^{\frac{1}{2}} \\
-X^{-\frac{3}{2}} Y^{{-\frac{1}{2}}} & +{4 X^{-\frac{1}{2}}}{Y^{{-\frac{1}{2}}}} &  -5 X^{{\frac{1}{2}}} Y^{{-\frac{1}{2}}}  &  +2 X^{{\frac{3}{2}}} Y^{{-\frac{1}{2}}} \\
 & & +X^{{\frac{1}{2}}} Y^{{-\frac{3}{2}}} &    - X^{{\frac{3}{2}}} Y^{{-\frac{3}{2}}}.
\end{array}
$$
The variables $X$ and $Y$ are the meridians to the oriented components
of the link as pictured in Figure~\ref{fig:9a42}.

\begin{figure}[ht]
\mbox{\vbox{\epsfbox{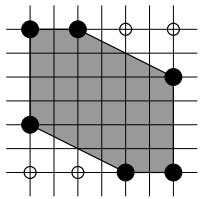}}}
\caption{\label{fig:9a42T} {\bf{Dual Thurston polytope $9a{42}$.}}
  The polytope is shaded. The four light circles are not in the dual
  Thurston polytope, according to the existence of the surface
  indicated in Figure~\ref{fig:9a42Seif}.}
\end{figure}

We claim that the dual Thurston polytope is the one pictured in
Figure~\ref{fig:9a42T}.  Of course, this is the Newton polytope of the
Alexander polynomial given above, and hence our claim is an immediate
consequence of Corollary~\ref{cor:AltLinks}.

However, the claim has also the following more elementary
proof. First, it is easy to find a disk which spans the component
$K_1$, meeting $K_2$ in four points. Puncturing the disk in these four points,
we obtain a surface in the link complement, from which we
conclude that the dual Thurston polytope is contained in the strip
$\{(x,y)\big| |x|\leq 3\}$.  Finding a similar disk spanning $K_2$,
we see that the Thurston polytope is contained in the square
$\{(x,y)\big| |x|\leq 3, |y|\leq 3\}$.

To narrow down the possibilities further, we use 
McMullen's bound, which states that the Thurston polytope contains the
Newton polytope of the multi-variable Alexander polynomial.  In view
of this, it remains to show that the dual Thurston polytope does not
contain any of the points $({3},{3})$,
$({1},{3})$, $(-{3},-{3})$, or
$(-{1},-{3})$.

\begin{figure}[ht]
\mbox{\vbox{\epsfbox{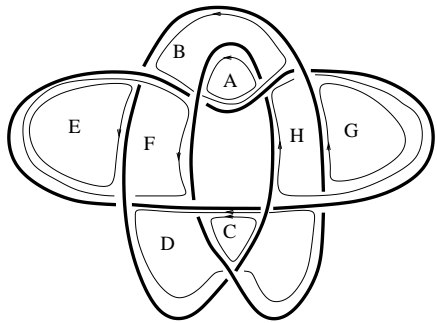}}}
\caption{\label{fig:9a42Seif}
{\bf{A representative of the homology class $(1,2)$.}}
This representative is obtained by first considering
the disks spanning the labeled circles, puncturing them in a minimal
number of points necessary, and then adding $9$ one-handles, to 
obtain the desired spanning surface.}
\end{figure}

But this follows at once from the fact that the homology class
$(1,2)$ is represented by a connected surface whose Euler
characteristic is $-5$. This surface is obtained as follows.
Consider the closed loops $A,...,H$ pictured in
Figure~\ref{fig:9a42Seif}.  These are to be thought of as closed loops
in the link complement.  Each bounds a disk in $S^3$; after puncturing
some of them -- namely, $B$, $H$, $F$, and $D$, each in a single point
-- we obtain parts of a surface. We attach one-handles to remove
the corner points from each of these disks, one at each crossing in
the projection (for example, $E$ and $B$ are connected at one
crossing; also, $C$ and $D$ are connected at one crossing). In this
manner, we obtain an immersed surface-with-boundary in the link
complement whose Euler characteristic is $-5$. This surface can be
readily resolved to obtain a smoothly embedded surface with the same
Euler characteristic, and representing the homology class $(1,2)$.

\subsection{A non-alternating example: $9^2_{50}=9n14$}

\begin{figure}[ht]
\mbox{\vbox{\epsfbox{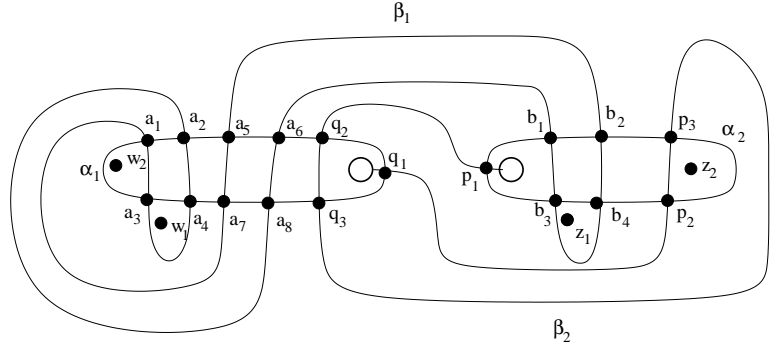}}}
\caption{\label{fig:9n14Diag}
{\bf{Heegaard diagram for $9n14$.}}}
\end{figure}

Consider the $9$-crossing link $9^2_{50}$ in Rolfsen's notation and
$9n14$ in Thistlethwaite's. This link is illustrated in
Figure~\ref{fig:9n14}. A $4$-pointed Heegaard diagram for this link can 
be drawn on a surface of genus one, as pictured in Figure~\ref{fig:9n14Diag}.
Inspecting this diagram, we see that the attaching circles $\alpha_i$ and
$\beta_j$ intersect in $36$ points as specified in the following table:

\begin{figure}[ht]
\mbox{\vbox{\epsfbox{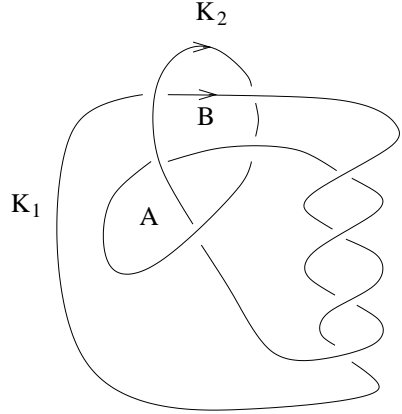}}}
\caption{\label{fig:9n14}
{\bf{The oriented link $9n14$.}}
Two regions in the projection complement, $A$ and $B$ are
distinguished.}
\end{figure}

\vskip.2cm
\begin{center}
\begin{tabular}{c|cc}
$\cap$ &  $\alpha_1$ & $\alpha_2$   \\
        \hline
$\beta_1 $ & $\{a_1,...,a_8\}$ & $\{b_1,...,b_4\}$ \\
$\beta_2$  & $\{q_1,q_2,q_3\}$ & $\{p_1,p_2,p_3\}$ 
\end{tabular}
\end{center}
\vskip.2cm

To calculate $\HH$-gradings, we use the following technique.  There is
a map
$$S^{i,j}\colon \alpha_i\cap\beta_j\longrightarrow \Z$$
defined as follows. $x,x'\in \alpha_i\cap\beta_j$, 
it is easy to see that there
are arcs $a\subset \alpha_i$ and $b\subset \gamma_i$, both going
from $x$ to $x'$, with the additional property that $a-b$ is homologous
to a sum of curves among the $\alpha_m$ and $\beta_n$. Let $D_{x,x'}$
be such a homology.
Then, $S^{i,j}$ is uniquely characterized
up to overall translation by the equation
$$S^{i,j}(x)-S^{i,j}(x')=(n_{z_1}(D_{x,x'})-n_{w_1}(D_{x,x'}),
n_{z_1}(D_{x,x'})-n_{w_1}(D_{x,x'})),$$
for all $x, x'\in\alpha_i\cap\beta_j$.
We say $\x,\y\in\Ta\cap\Tb$ have the same {\em type}, 
if  there is some reordering $\sigma\in S_2$ (the symmetric
group on two letters) so that $\x=(x_1,x_2)$ and $\y=(y_1,y_2)$ 
and $x_i,y_i\in \alpha_i\cap\beta_{\sigma(i)}$ for $i=1,2$.
Now suppose that $\x,\y\in\Ta\cap\Tb$ have the same type, and $\sigma$
is the corresponding transposition, then
according to Equation~\eqref{eq:DefHGrading},
$$\HGrading_{\w,\z}(x_1\times x_2)-\HGrading_{\w,\z}(y_1\times y_2)
=S^{1,\sigma(1)}(x_1)+S^{2,\sigma(2)}(x_2)-S^{1,\sigma(1)}(y_1)
-S^{2,\sigma(2)}(y_2).$$
Thus, to determine $\HGrading_{\w,\z}$ up to overall translation,
it suffices to calculate all the $S^{i,j}$, and then connect up two
intersection points of different types. We call $S^{i,j}(x)-S^{i,j}(x')$
the {\em relative difference} between $x$ and $x'$, and drop the 
superscript $i,j$ from the notation.

\begin{figure}[ht]
\mbox{\vbox{\epsfbox{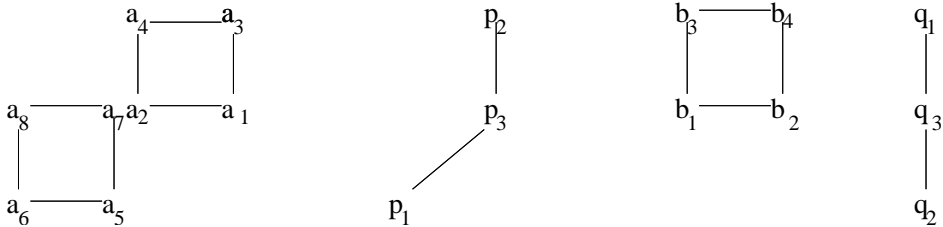}}}
\caption{\label{fig:9n14RelDiff}
{\bf{Relative differences of intersection points for $9n14$.}}
We plot the
relative differences of the intersection points of types $a_i$, $p_j$,
$b_k$, and $q_\ell$. Edges represent two comparable intersection points; for
example, $S(a_3)-S(a_4)=(1,0)$, $S(p_3)-S(p_1)=(1,1)$.}
\end{figure}

\begin{figure}[ht]
\mbox{\vbox{\epsfbox{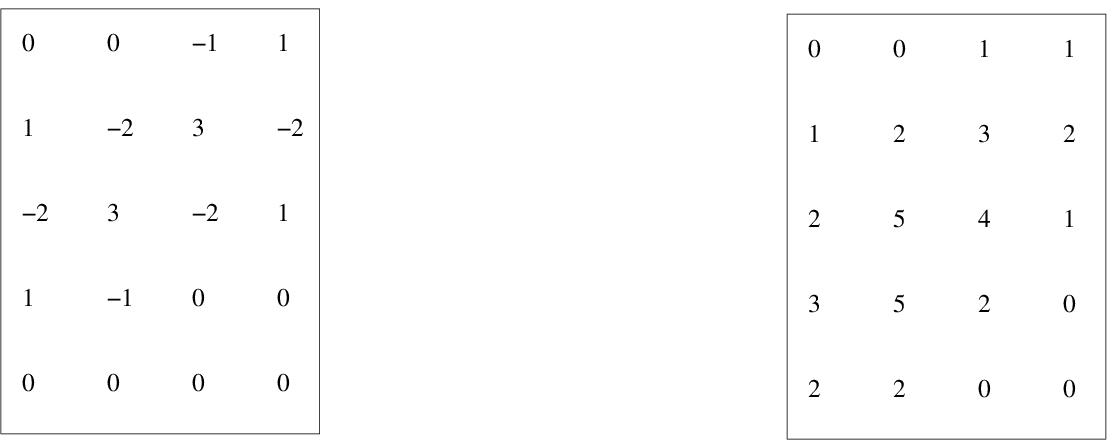}}}
\caption{\label{fig:9n14Ans}
{\bf{Ranks and Euler characteristics of $\CFLa$ for $9n14$.}}}
\end{figure}

Relative differences of intersection points $\alpha_i\cap\beta_j$
can also be found, as in Figure~\ref{fig:9n14RelDiff}. Moreover,
it is easy to find a small square in Figure~\ref{fig:9n14Diag}
which contains none of the basepoints, connecting
$a_6\times p_1$ and $b_1\times q_2$. With this information, now,
it is straightforward to calculate the $\HH$-gradings of all 
the generators of $\Ta\cap\Tb$. In fact, in Figure~\ref{fig:9n14Ans}
we have displayed the ranks of $\CFLa$ in each different $\HH$-grading,
and next to it the Euler characteristics of these groups.

Using this information, together with Equation~\eqref{eq:Symmetry}, we
conclude immediately that for each $h\in \HH$, the rank of
$\HFLa(9n14,h)$ coincides with the absolute value of its Euler
characteristic. In particular, thanks to
Theorem~\ref{thm:ThurstonNorm}, the Newton polygon of the Alexander
polynomial  and the dual Thurston polytope of $9n14$ coincide.

Note that an alternative argument to determine the Thurston polytope
can be given as follows. First, observe that both components $K_1$ and
$K_2$ can be spanned by the surfaces in the link complement which
have $\chi=-2$ (in the case of $K_1$, we use a disk with three
punctures, in the case of $K_2$, we use a torus with two
punctures). The fact, now, that the Thurston polytope is no larger
than the Newton polytope follows from the fact that the relative
homology class $(1,-1)$ is Poincar\'e dual to a fibration. This can be
seen using Gabai's method of disk
decompositions~\cite{DiskDecomposition}. Specifically, consider the
Seifert surface $F$ for $9n14$ obtained as follows. Consider the
checkerboard coloring of the link projection where the regions $A$ and
$B$ are colored white. The black regions can be used to construct the
Seifert surface $F$, and consider the corresponding sutured
manifold. Now, attaching a disk along $A$, and then one along $B$
(which meet the sutures in two points apiece), we end up with the
sutured manifold consisting of the solid torus with two parallel
sutures which are meridians. Attaching one more disk $C$, we end up
with a three-ball with a single suture along the equator. Since each
of our disks $A$, $B$, and $C$ met the sutures along two points
apiece, and since we end up with the three-ball with a single suture,
it follows that we started with a surface $F$ which is the fiber of a
fibration of the link complement.

\subsection{The link $9^3_{21}=9n27$}
\begin{figure}[ht]
\mbox{\vbox{\epsfbox{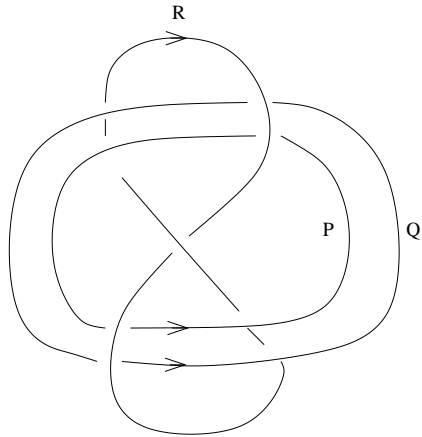}}}
\caption{\label{fig:9n27} {\bf{The link $9n27$.}}}
\end{figure}

Consider the link $9n27$ considered in Figure~\ref{fig:9n27}, the
nine-crossing, non-split link with trivial multi-variable Alexander polynomial.
We can draw a compatible Heegaard diagram on the sphere, as illustrated
in Figure~\ref{fig:9n27}.
\begin{figure}[ht]
\mbox{\vbox{\epsfbox{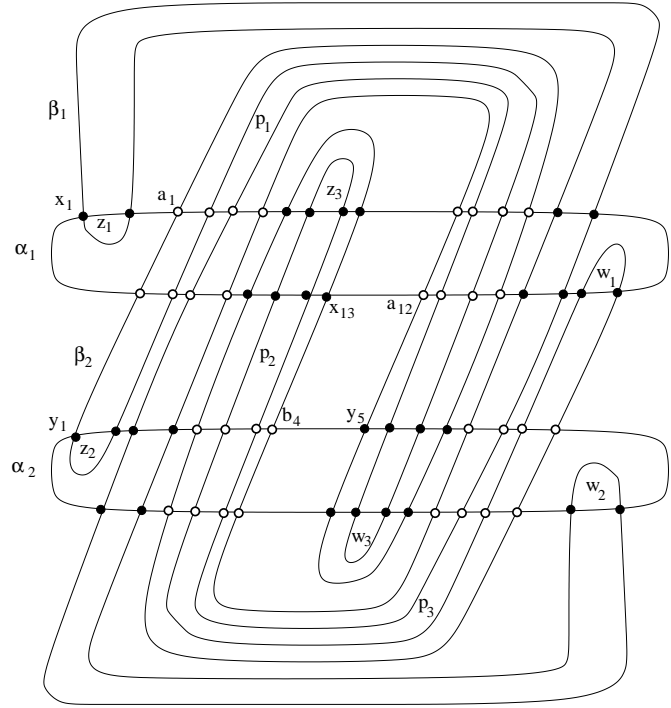}}}
\caption{\label{fig:9n27Diag} {\bf{Heegaard diagram for $9n27$.}}
  The intersection points $a_i$, $b_i$, $x_i$, and $y_i$ are ordered
  in a clockwise order as we traverse the corresponding $\alpha$
  circle.  (The additional basepoints points $p_i$ $i=1,2,3$ will be
  used in the calculation of the link Floer homology groups, but
  should be ignored for the moment.)}
\end{figure}

\begin{figure}[ht]
\mbox{\vbox{\epsfbox{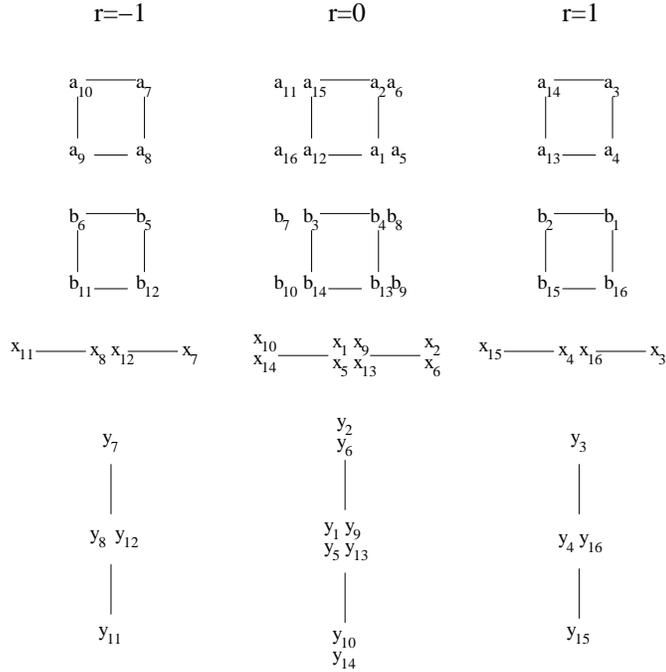}}}
\caption{\label{fig:9n27Diffs} {\bf{Differences of intersection points
for $9n27$.}}
The horizontal component measures the $p$-coordinate, the
vertical the $q$ coordinate, and the three different columns
correspond to different $r$ coordinates, as indicated.}
\end{figure}

This has two pairs of  attaching circles $\{\alpha_1,\alpha_2\}$
and $\{\beta_1,\beta_2\}$, which intersect according to the pattern
illustrated in the following table
\vskip.2cm
\begin{center}
\begin{tabular}{c|cc}
$\cap$ &  $\alpha_1$ & $\alpha_2$   \\
        \hline
$\beta_1 $ & $\{x_1,...,x_{16}\}$ & $\{b_1,...,b_{16}\}$\\
$\beta_2$  &  $\{a_1,...,a_{16}\}$ & $\{y_1,...,y_{16}\}$  \\
\end{tabular}
\end{center}
\vskip.2cm
Thus, there are $512$ intersection points of $\Ta\cap\Tb$,
of the two types $a_i\times b_j$ and $x_i\times y_j$ with
$i,j\in\{1,...,16\}$. Relative differences between
intersection points are given as in Figure~\ref{fig:9n27Diffs},
with the convention that if $\x$ and $\y$ are two intersection points,
then 
$$(p,q,r)=(n_{z_1}(\phi)-n_{w_1}(\phi),
n_{z_2}(\phi)-n_{w_2}(\phi),
n_{z_3}(\phi)-n_{w_3}(\phi))$$
for $\phi\in\pi_2(\x,\y)$. (The components of $9n27$ corresponding to
$z_1$, $z_2$, and $z_3$ are labelled by $P$, $Q$, and $R$ respectively
in Figure~\ref{fig:9n27}.)

\begin{figure}[ht]
\mbox{\vbox{\epsfbox{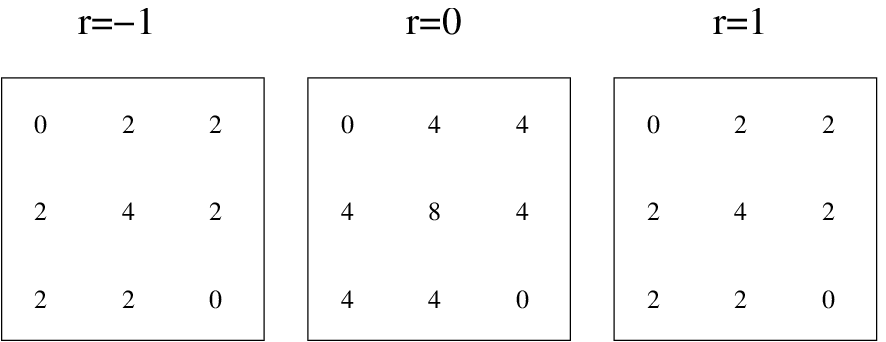}}}
\caption{\label{fig:9n27Ans} {\bf{Ranks of $\HFLa$ for $9n27$.}}}
\end{figure}

We aim to show that the ranks of the various groups are 
given in Figure~\ref{fig:9n27Ans}. Indeed, for the purposes
of calculating the link Floer homology polytope, it suffices
to verify this calculation in the cases where
$(p,q,r)\neq (0,0,0)$.

To verify these calculations, observe that there is a collection
of obvious small rectangles, giving flowlines pairing off
\[
\begin{array}{ll}
x_i\times y_j ~\text{and}~ a_j\times b_{i-2}
&{\text{if $i=3,...,8$ and $j=1,...,8$}} \\
x_i \times y_j ~\text{and}~ a_{j-2}\times b_i
& {\text{if $i=9,...,16$ and $j=11,...,16$}} \\
x_i\times y_j ~\text{and}~ a_{17-j}\times b_{17-i} 
& {\text{if $i=9,...,16$~\text{and}~$j=1,...,8$}} \\
x_i\times y_j ~\text{and}~ a_{19-j}\times b_{19-i} 
& {\text{if $i=3,...,8$~\text{and}~$j=11,...,16$}} \\
\end{array}
\]
For example, there is a rectangle
giving rise to a pseudo-holomorphic disk
connecting $a_{12}\times b_4$ to $x_{13}\times y_5$; also, there
is a similar rectangle going from $x_{16}\times y_4$ to $a_{13}\times b_{1}$.

These flows prove that for some $\HH$-grading $(p,q,r)$ for
which there are no generators containing $x_1$, $x_2$, 
$y_9$, or $y_{10}$, then the corresponding homology group vanishes
$\HFLa(L,(p,q,r))=0$.

In particular, it follows from this observation, together with the
calculation of relative differences from Figure~\ref{fig:9n27Diffs} that
\begin{eqnarray}
\HFLa(L,(p,q,r))=0
& {\text{if $|p|$, $|q|$, or $|r|$  is greater than $1$}} \nonumber \\
& {\text{or if $(p,q)\in\{(1,1),(-1,-1)\}$}}. \label{eq:VanishingVertices}
\end{eqnarray}

With some additional work, we now show that 
\begin{equation}
\label{eq:OneVertex}
\HFLa(L,(-1,-1,-1))=2.
\end{equation} 
This calculation is elementary, using only the property that
any rectangle has a pseudo-holomorphic representative, and also that 
homotopy classes which have negative multiplicity somewhere admit no
pseudo-holomorphic representatives. We organize this as follows.

There are eight generators of $\CFLa(L,(-1,-1,-1))$,
$$\left\{
\begin{array}{llll}
x_{10}\times y_{11}, & x_{11}\times y_{10}, &x_{11}\times y_{14}, & x_{14}\times y_{11}, \\
a_9\times b_{14},& a_{9}\times b_{10},& a_{12}\times b_{11},& 
a_{16}\times b_{11}
\end{array}\right\}.$$
It is straightforward to find six rectangles connecting such generators
whose local multiplicities
vanish at $\{w_1,w_2,w_3,z_1,z_2,z_3\}$. Indeed, we assemble our generators into
three sets, under an equivalence relation generated by the rectangles,
\begin{eqnarray*}
A&=&\{a_{16}\times b_{11}\}, \\
B&=&\{x_{10}\times y_{11}, x_{11}\times y_{10}, a_9\times b_{10}\} \\
C&=&\{x_{11}\times y_{14},a_9\times b_{14}, a_{12}\times b_{11},
x_{14}\times y_{11}\}.
\end{eqnarray*}
Observe that any homotopy class connecting an element of $C$ to an
element of $A$ or $B$ must have positive local multiplicity somewhere.
(To this end, it is useful to make the following observations.
Consider the pointS $p_1$ and $p_3$ in Figure~\ref{fig:9n27Diag}.
Observe that it has multiplicity zero in any periodic domain, and also
in any of the six rectangles pictured above. Thus, it suffices to find
two homotopy class connecting some intersection point of $C$ to one
$A$ and $B$ respectively, and to verify that for each, the sum of the
local multiplicities at $p_1$ and $p_3$ are positive. This is
straightforward.)  It follows that $A$ and $B$ generate a subcomplex
of $\CFLa(L,(-1,-1,-1))$, with quotient complex $C$. The complex $C$
is easily seen to have trivial homology (four of the six rectangles
connect generators in $C$ in such a manner that it has trivial
homology). Moreover, it is easy to see that any homotopy class
connecting a generator of $A$ with a generator of $B$ must have both
positive and negative local multiplicities (using both basepoints
$p_2$ and $p_3$ from Figure~\ref{fig:9n27Diag}).  Obviously $H_*(A)$
has rank one. Moreover, the two rectangles connecting generators of
$B$ show that the boundary operator on this latter complex is given by
$$\partial(x_{10}\times y_{11})=a_9\times b_{10}=\partial x_{11}\times
y_{10},$$ so that $H_*(B)$ has rank one, as well.

We could proceed in this manner to analyze the other $(p,q,r)$-levels.
There is, however, a quicker argument which allows us to determine the
link Floer homology polytope with little extra work (proving
non-triviality of
$\HFLa(L,(p,q,r))$ at the remaining
points $(p,q,r)$
where at most one of the coordinates
vanishes, building on the calculations from
Equations~\eqref{eq:VanishingVertices} and \eqref{eq:OneVertex}.

If we consider an isotopic translate $\beta_2'$ of $\beta_2$, where we
allow the isotopy to cross $z_1$.  The remaining homology retains its
$(q,r)$-grading.  We can arrange for this isotopy to eliminate all
intersections containing $\{a_i\}_{i=1}^{16}$, leaving only $x_{9}$ and
$x_{10}$. Thus, in grading $(q,r)=(-1,-1)$, there is a Floer homology
group $G$ generated by the two remaining intersection points $[x_9\times
y_{11}]$ and $[x_{10}\times y_{11}]$. 

By isotopy invariance of homology, this result can be interpreted as
follows. Consider the enhanced differential $D_1\colon
\CFLa(L)\longrightarrow \CFLa(L)$ defined as in
Equation~\eqref{eq:DefD}, only now allowing disks with
$n_{z_1}(\phi)\neq 0$. This complex has a remaining $(q,r)$ grading, 
and it is isomorphic to the Floer homology theory defined using
$\beta_2'$ in place of $\beta_2$. 
(Although this might appear to be an {\em ad hoc} construction,
this is in fact an example of a more general principle investigated in
\cite[Section~\ref{Links:subsec:Forgetfuls}]{Links}: the homology
groups of a differential counting
holomorphic disks which cross a basepoint $z_i$ corresponds to
the homology groups of the link obtained by removing the $i^{th}$ component,
and then tensoring with a vector space of rank two.)
In particular, for the grading level $(q,r)=(-1,-1)$,
we obtain the long exact sequence
$$
\begin{CD}
...@>>> \HFLa(L(-1,-1,-1))@>>> G@>>> \HFLa(L(0,-1,-1)) @>>>...
\end{CD}
$$
Now, both $\HFLa(L(-1,-1,-1))$ and $G$ are two dimensional, but
their generators have different Maslov gradings. It follows at once
that $\HFLa(L(0,-1,-1))$ must be non-trivial.

Combining the above with Equation~\eqref{eq:OneVertex},
and then using Equation~\eqref{eq:Symmetry} and the additional symmetries
\begin{align*}
\HFLa(L,(p,q,r))&\cong \HFLa(L,(q,p,r)) \\
\HFLa(L,(p,q,r))&\cong \HFLa(L,(p,q,-r))
\end{align*}
 (which follow from
corresponding symmetries of the link), it is now straightforward to
verify that $\HFLa(L,(p,q,r))\neq 0$ for 
$$(p,q,r)\in \left\{
\begin{array}{llllll}
(-1,-1,-1),&
(0,-1,-1), &(-1,0,-1),&(1,1,-1),&(0,1,-1),&(1,0,-1), \\
(1,1,1),&(0,1,1),&(1,0,1),&
(-1,-1,1),&(0,-1,1),&(-1,0,1)
\end{array}
\right\},
$$ giving the link Floer homology polytope. (A more tedious calculation
along similar lines can be used to verify that the ranks of the Floer homology
groups of $9n27$ are as given in Figure~\ref{fig:9n27Ans}.)

\begin{figure}[ht]
\mbox{\vbox{\epsfbox{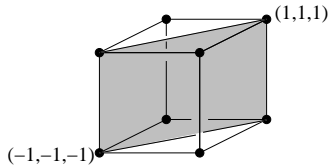}}}
\caption{\label{fig:9n27T} {\bf{Thurston polytope for $9n27$.}}
The Thurston polytope is the lightly shaded region.}
\end{figure}

This result is consistent with Thurston polytope pictured in
Figure~\ref{fig:9n27T}. (For $9n27$, it is easy to give
a more direct verification of its 
Thurston polytope.)

\subsection{Kinoshita-Terasaka links}

\begin{figure}[ht]
\mbox{\vbox{\epsfbox{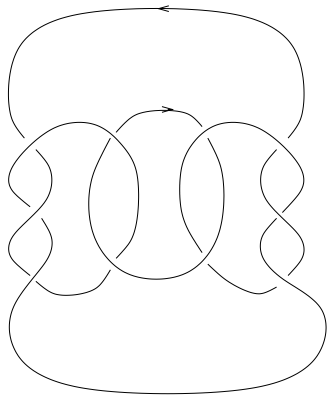}}}
\caption{\label{fig:KTLink}
{\bf{The Kinoshita-Terasaka link.}}  This link (denoted $L10n36$
in Thistlethwaite's notation) has vanishing Alexander polynomial. }
\end{figure}

\begin{figure}[ht]
\mbox{\vbox{\epsfbox{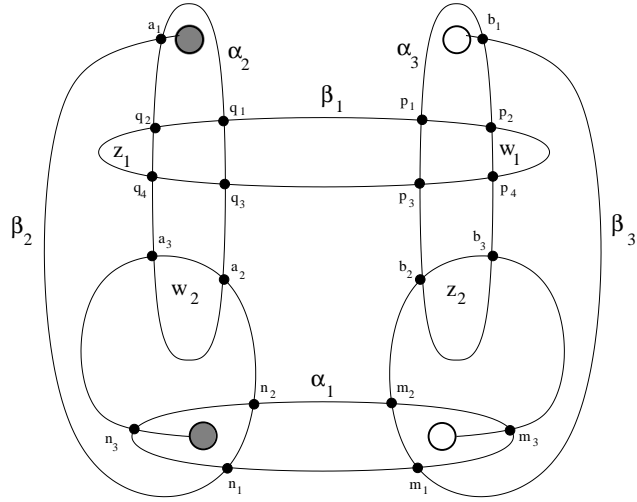}}}
\caption{\label{fig:KTDiag}
{\bf{Heegaard diagram for the Kinoshita-Terasaka link.}} A
four-pointed Heegaard diagram for the link. This picture
takes place on the surface of genus two, obtained by identifying
the two white circles and the two gray circles, and adding a point
at infinity. The pair of basepoints
$w_1$ and $z_1$ represent the unknotted component, and $w_2$ and $z_2$
represent the connected sum of two trefoils.}
\end{figure}

Consider the Kinoshita-Terasaka link pictured in
Figure~\ref{fig:KTLink}.  This link has the $4$-pointed Heegaard
diagram pictured in Figure~\ref{fig:KTDiag}. The Heegaard surface has
genus $2$, and is equipped with attaching circles
$\{\alpha_1,..,\alpha_3\}$ and $\{\beta_1,...,\beta_3\}$.  These
attaching circles meet in intersection points according to the
following table.

\vskip.2cm
\begin{center}
\begin{tabular}{c|ccc}
$\cap$ &  $\alpha_1$ & $\alpha_2$ & $\alpha_3$  \\
        \hline
$\beta_1 $ & $\emptyset$ & $\{q_1,...,q_4\}$ & $\{p_1,...,p_4\}$\\
$\beta_2$  & $\{n_1,...,n_3\}$ & $\{a_1,...,a_3\}$ & $\emptyset$ \\
$\beta_3$ & $\{m_1,...,m_3\}$ & $\emptyset$ & $\{b_1,...,b_3\}$
\end{tabular}
\end{center}
\vskip.2cm

Thus, intersection points $\Ta\cap\Tb$ have the form
$a_i\times m_j \times p_k$ and $b_i\times n_j \times q_k$,
with $i=1,...,3$, $j=1,...,3$ and $k=1,...,4$.

To calculate the relative values of $\HGrading_{\w,\z}$, we
first calculate the relative differences
$(n_{z_1}(\phi)-n_{w_1}(\phi), n_{z_2}(\phi)-n_{w_2}(\phi))$,
where the homotopy classes $\phi$ are Whitney disks between
the various intersection points of the $\alpha_i$ and the $\beta_j$.

\begin{figure}[ht]
\mbox{\vbox{\epsfbox{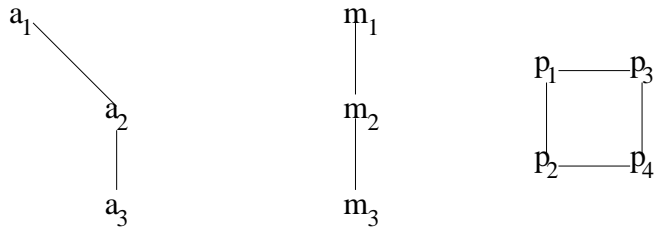}}}
\caption{\label{fig:RelDiff}
{\bf{Relative differences of intersection points.}}  
The relative differences for $b_i$, $n_j$ and $q_k$ have the same shape,
replacing the symbols $a$, $m$, and $p$ by $b$, $n$, and $q$ respectively.}
\end{figure}

These relative differences, for intersection points of the form $a_i$,
$m_j$, and $p_k$ are plotted in Figure~\ref{fig:RelDiff}.  Relative
differences of the $b_i$, $n_j$, and $q_k$ have the same corresponding
shapes.

Indeed, we now claim that $$\HGrading(a_i\times m_j\times p_k)
=\HGrading(b_i\times n_j\times p_k)$$
for all $i=1,...,3$, $j=1,...,3$, and $k=1,...,4$. This is exhibited, for example, by the
obvious ``large hexagon'' containing the point at infinity in
Figure~\ref{fig:KTLink}, thought of as a Whitney disk from $a_1\times
m_1\times p_1$ to $b_1\times n_1\times q_1$ (which is, of course,
disjoint from $\ws$ and $\zs$.

\begin{figure}[ht]
\mbox{\vbox{\epsfbox{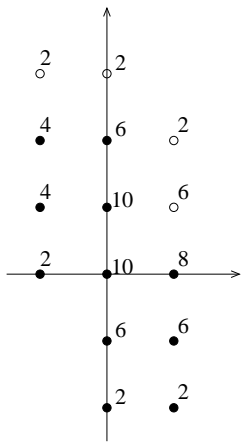}}}
\caption{\label{fig:RankChainGroups}
{\bf{Ranks of chain groups for the Kinoshita-Terasaka link.}}  We plot the
rank of the chain complex $\CFLa(\orL,(i,j))$ for
each $(i,j)$ coming from the above diagram.
The upper left-hand-corner is generated by $a_1\times m_1\times p_1$
and $b_1\times n_1\times q_1$. The three empty dots represent 
levels where, although this chain complex is non-trivial,
the homology $\HFLa$ is trivial.}
\end{figure}

Thus, the ranks of the non-zero chain groups in each (relative)
$\Z^2$-grading are as illustrated in Figure~\ref{fig:RankChainGroups}.

\begin{figure}[ht]
\mbox{\vbox{\epsfbox{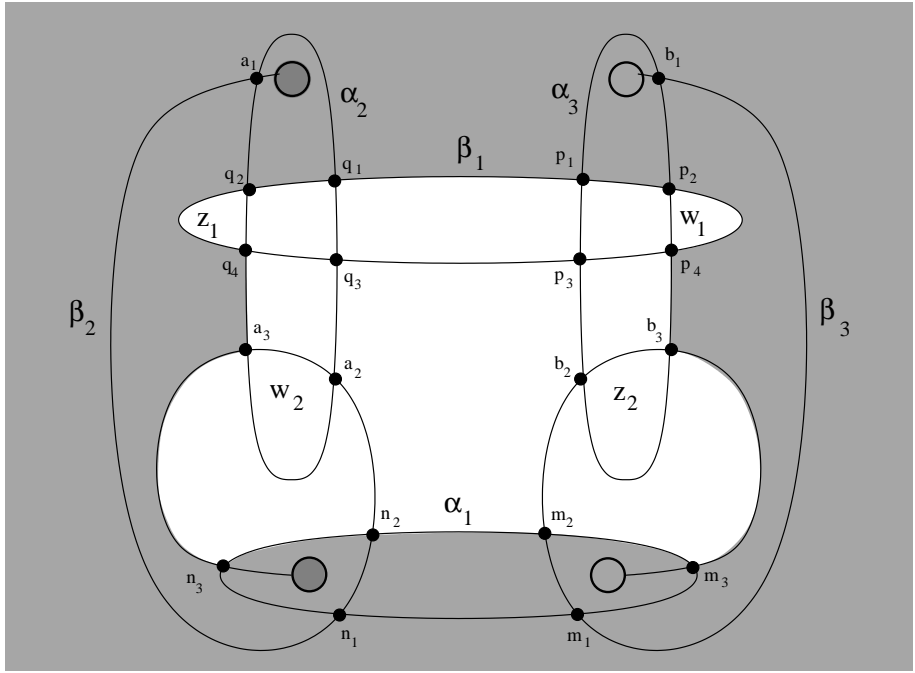}}}
\caption{\label{fig:TheDomain}
{\bf{A domain connecting $a_3\times m_3\times p_4$ and $b_3\times
n_3\times q_4$.}}  The only two-chain with local multiplicties $0$ or $+1$
(indicated by gray shading) representing 
a homotopy class $\phi\in\pi_2(b_3\times n_3\times q_4,a_3\times m_3\times
p_4)$.}
\end{figure}

We claim that in the filtration level given by $(1,-2)$, the two
generators $a_3\times m_3\times p_4$ and $b_3\times n_3\times q_4$
both survive in homology. This is true because the Maslov grading
of $a_3\times m_3\times p_4$ is one greater than that of
$b_3\times n_3\times q_4$, and there are no positive domains from
the first generator to the second.
This can be seen at once by considering the
domain connecting them illustrated in Figure~\ref{fig:TheDomain}.
This is a domain from 
More specifically, that domain gives a domain
$\phi\in\pi_2(b_3\times n_3\times q_4,
a_3\times m_3\times p_4)$ with no negative local multiplicities
(and Maslov index $-1$). From this, it is easy to see that there are
no nowhere negative domains which miss all basepoints and
have with Maslov index one between the two generators, and hence
no differentials.

\begin{figure}[ht]
\mbox{\vbox{\epsfbox{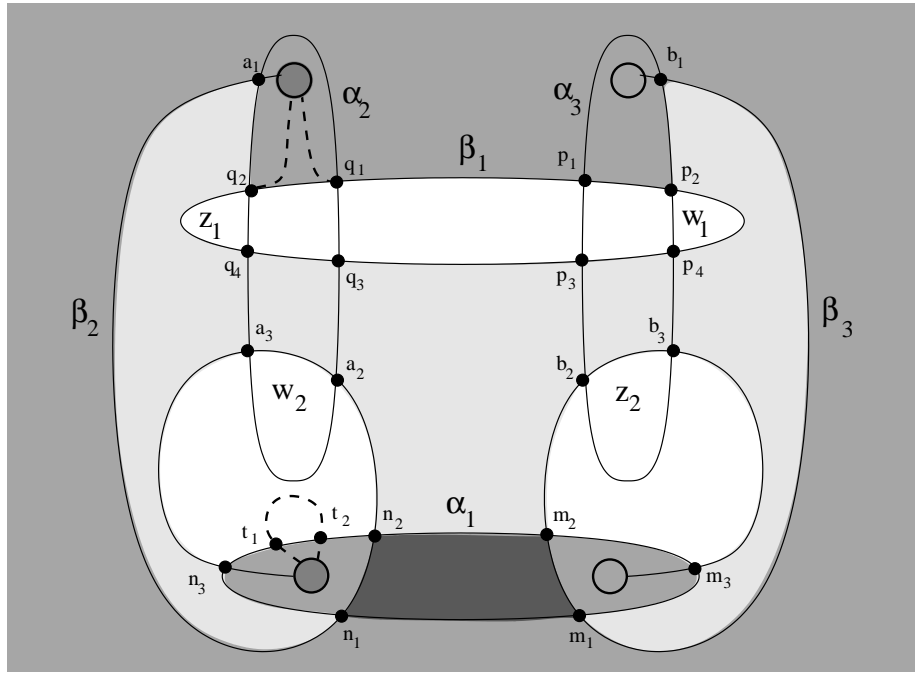}}}
\caption{\label{fig:Finger} {\bf{A domain connecting $a_1\times
      m_3\times p_2$ and $b_1\times n_3\times q_2$.}}  This domain has
  Maslov index equal to $+1$, and positive local
  multiplicities. (Here, the multiplicities run between $0$ and $3$:
  darker shading means higher local multiplicity.) Performing the
  finger move indicated by the dotted arc, we obtain an isotopic copy
  of $\beta_1$ which meets $\alpha_1$ in two points $t_1$ and
  $t_2$. After this finger move, the resulting domain acquires some
  negative local
  multiplicity $-2$.}
\end{figure}

We claim also that in the filtration level given by $(-1,0)$, the two
generators $a_1\times m_3\times p_2$ and $b_1\times n_3\times q_2$
survive in homology. This is slightly more subtle, since there now is
a (unique) domain connecting $a_1\times m_3\times p_2$ to $b_1\times
n_3\times q_2$, which has only non-negative local multiplicities (and
Maslov index one). This domain is pictured in
Figure~\ref{fig:Finger}. Perform an isotopy of the diagram, the
``finger move'' indicated by the dotted line in the figure, which
introduces two new intersection points $t_1$ and $t_2$ of $\alpha_1$
with $\alpha_2$. After this isotopy, all domains connecting the two
fixed generators have both positive and negative multiplicities. It is
important to note, though, that the isotopy introduces $18$ new
generators, of the form $t_i\times a_j \times b_k$ $i\in\{1,2\}$,
$j,k\in\{1,2,3\}$. However, using the small rectangle supported in the
finger connecting $a_1\times t_1$ to $n_3\times q_2$ (which crosses no
basepoints), we see that $t_i\times a_1\times a_2$ is in the same
bigrading as $b_1\times n_3\times q_2$, and hence that none of the
newly-introduced intersection points is supported in $\HH$-grading
$(-1,0)$.

With this input, together with the usual symmetry property
(Equation~\eqref{eq:Symmetry}), the link Floer homology polytope is
determined immediately. In particular, it follows that the relative
filtration levels displayed in Figure~\ref{fig:RankChainGroups}
indeed coincide with the absolute $\HH$-grading, and hence also that
the homology groups in $\HH$-grading $(-1,3)$, $(0,3)$, $(1,2)$, and
$(1,1)$ are in fact trivial.

\begin{figure}[ht]
\mbox{\vbox{\epsfbox{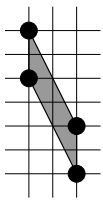}}}
\caption{\label{fig:KTT}
{\bf{Dual Thurston polytope of the Kinoshita-Terasaka link.}}  This is
the dual Thurston polytope of the link pictured in
Figure~\ref{fig:KTLink}, where the horizontal axis is represented by
multiples of the meridian of the unknot component, and the vertical
axis is represented by multiples of the meridian for the component
which is a connected sum of trefoils.  These meridians inherit
orientations as indicated in Figure~\ref{fig:KTLink}.}
\end{figure}

We conclude from this, together with Theorem~\ref{thm:ThurstonNorm}
that the dual Thurston polytope for the Kinoshita-Terasaka link is as
pictured in Figure~\ref{fig:KTT}. In particular, this suggests that
the homology class dual to $(1,-2)$ is represented by a surface $F$
with $\chi(F)=-1$. It is now an exercise in visualization to find such
a representative (given by a sphere with three punctures).

\begin{figure}[ht]
\mbox{\vbox{\epsfbox{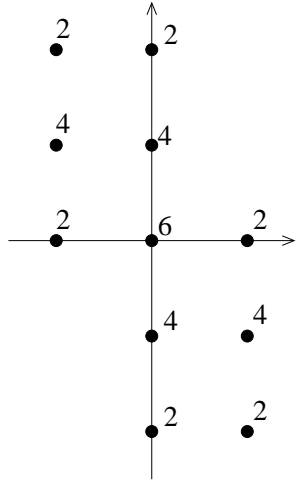}}}
\caption{\label{fig:KTAns}
{\bf{Ranks of link Floer homology  groups
 for the Kinoshita-Terasaka link.}}  We plot the
rank of the groups $\HFLa(\orL,(i,j))$ for
each $(i,j)$.}
\end{figure}

Again, a more involved calculation using the same circle of ideas can
be used to calculate the full link Floer homology groups of the
Kinoshita-Terasaka link, given in Figure~\ref{fig:KTAns}.

\bibliographystyle{plain}
\bibliography{biblio}

\end{document}